\pdfoutput=1
\RequirePackage{ifpdf}
\ifpdf %We are running pdfTeX in pdf mode
\documentclass[pdftex]{sigma}
\else
\documentclass{sigma}
\fi

\usepackage[all]{xy}
\usepackage{epic,eepic,multicol}
\numberwithin{equation}{section}
\newtheorem{thm}{Theorem}[section]
\newtheorem{prop}[thm]{Proposition}
\newtheorem{lem}[thm]{Lemma}

{\theoremstyle{definition}
\newtheorem{defn}[thm]{Definition}
\newtheorem{rem}[thm]{Remark}
\newtheorem{ex}[thm]{Example}}

\def\topseudosquare(#1,#2,#3){
\hbox{$\hskip-0.4pt\vtop to #1{\normalbaselines\m@th
\hrule\vfil\hbox to #2{$\hfill\scriptstyle #3\hfill$}\vfil\hrule}$\hskip-0.4pt
\vrule}}

\def\Topseudosquare<#1,#2,#3>{
\hbox{$\hskip-0.4pt\vtop to #1{\normalbaselines\m@th
\hrule\vfil\hbox to #2{$\hfill\scriptstyle #3\hfill$}\vfil\hrule}$\hskip-0.4pt
\vrule}}

\def\TTopseudosquare<#1,#2,#3>{
\hbox{$\hskip-0.4pt\vtop to #1{\normalbaselines\m@th
\hrule\vfil\hbox to #2{$\hfill #3\hfill$}\vfil\hrule}$\hskip-0.4pt
\vrule}}

\def\tophalfsquare(#1,#2,#3){
\hbox{\vrule$\hskip-0.4pt\vtop to #1{\normalbaselines\m@th
\hrule\vfil\hbox to #2{\hfill$\scriptstyle #3$\hfill}\vfil\hrule}$\hskip-1pt
\vrule}}

\def\Tophalfsquare<#1,#2,#3>{
\hbox{\vrule$\hskip-0.4pt\vtop to #1{\normalbaselines\m@th
\hrule\vfil\hbox to #2{\hfill$\scriptstyle #3$\hfill}\vfil\hrule}$\hskip-1pt
\vrule}}

\def\TTophalfsquare<#1,#2,#3>{
\hbox{\vrule$\hskip-0.4pt\vtop to #1{\normalbaselines\m@th
\hrule\vfil\hbox to #2{\hfill$#3$\hfill}\vfil\hrule}$\hskip-1pt
\vrule}}

\def\TY(#1,#2,#3){#1^{(#2)}_{#3}}

\def\m@th{\mathsurround=0pt}

\DeclareMathOperator{\sgn}{sgn}
\DeclareMathOperator{\diag}{diag}
\begin{document}

\newcommand{\arXivNumber}{1409.8622}

\allowdisplaybreaks

\renewcommand{\thefootnote}{$\star$}

\renewcommand{\PaperNumber}{033}

\FirstPageHeading

\ShortArticleName{Cluster Variables and Monomial Realizations of Crystal Bases}

\ArticleName{Cluster Variables on Certain Double Bruhat Cells\\
of Type $\boldsymbol{(u,e)}$ and Monomial Realizations\\
of Crystal Bases of Type~A\footnote{This paper is a~contribution to the Special Issue on New Directions in Lie Theory.
The full collection is available at \href{http://www.emis.de/journals/SIGMA/LieTheory2014.html}{http://www.emis.de/journals/SIGMA/LieTheory2014.html}}}

\Author{Yuki KANAKUBO and Toshiki NAKASHIMA}

\AuthorNameForHeading{Y.~Kanakubo and T.~Nakashima}

\Address{Division of Mathematics, Sophia University, Yonban-cho 4, Chiyoda-ku,\\
Tokyo 102-0081, Japan}
\Email{\href{mailto:j_chi_sen_you_ky@sophia.ac.jp}{j\_chi\_sen\_you\_ky@sophia.ac.jp}, \href{mailto:toshiki@sophia.ac.jp}{toshiki@sophia.ac.jp}}

\ArticleDates{Received October 01, 2014, in f\/inal form April 14, 2015; Published online April 23, 2015}

\Abstract{Let~$G$ be a~simply connected simple algebraic group over $\mathbb{C}$,~$B$ and $B_-$ be two opposite Borel
subgroups in~$G$ and~$W$ be the Weyl group.
For~$u$, $v\in W$, it is known that the coordinate ring ${\mathbb C}[G^{u,v}]$ of the double Bruhat cell
$G^{u,v}=BuB\cap B_-vB_-$ is isomorphic to an upper cluster algebra $\bar{{\mathcal A}}({\bf i})_{{\mathbb C}}$ and
the generalized minors $\{\Delta(k;{\bf i})\}$ are the cluster variables belonging to a~given initial seed in
${\mathbb C}[G^{u,v}]$ [Berenstein~A., Fomin~S., Zelevinsky~A., \textit{Duke Math.~J.} \textbf{126} (2005), 1--52].
In the case $G={\rm SL}_{r+1}({\mathbb C})$, $v=e$ and some special $u\in W$, we shall describe the generalized minors
$\{\Delta(k;{\bf i})\}$ as summations of monomial realizations of certain Demazure crystals.}


\Classification{13F60; 81R50; 17B37}

\renewcommand{\thefootnote}{\arabic{footnote}}
\setcounter{footnote}{0}

\section{Introduction}

As is well-known that theory of cluster algebras has been initiated by S.~Fomin and A.~Zelevinsky in the study of product
expressions by~$q$-commuting elements for upper global bases (=~dual canonical bases).
Crystal bases are obtained from global bases considering the parameter~$q$ at~0.
Thus, we can guess that they should be deeply related each other at their origins.

Let~$G$ be a~simply connected simple algebraic group over $\mathbb C$ of rank~$r$.
Let~$B$ and $B_-$ be two opposite Borel subgroups in~$G$, $N\subset B$ and $N_-\subset B_-$ their unipotent radicals,
$H:=B\cap B_-$ a~maximal torus, and~$W$ the associated Weyl group.
In~\cite{A-F-Z}, it is shown that for $u,v\in W$ the coordinate ring $\mathbb C[G^{u,v}]$ of double Bruhat cell
$G^{u,v}:=BuB\cap B_-vB_-$ has the structure of an upper cluster algebra.
The initial cluster variables of this upper cluster algebras are given as certain generalized minors on $G^{u,v}$.

In~\cite{N1}, the second author revealed the relations between some generalized minors and monomial realizations of
crystal bases.
A~naive def\/inition of monomial realizations of crystal bases is as follows (see Section~\ref{section3} for the exact def\/initions): Let
${\mathcal Y}$ be the set of monomials in inf\/initely many variables (see Section~\ref{section3}, equation~\eqref{cydef}).
We shall def\/ine the crystal structures on ${\mathcal Y}$ associated with certain set of integers $p=(p_{i,j})_{1\leq
i\ne j\leq r}$ and a~Cartan matrix.
And we can obtain a~crystal for an irreducible module as a~connected component of ${\mathcal Y}$.
For example, for type $A_4$ and $p_{i,j}=1$ if $i<j$ and $p_{i,j}=0$ if $i>j$, we have the following crystal graph of
the crystal~$B(\Lambda_3)$, where $\Lambda_3$ is the 3$^{\rm rd}$ fundamental weight.
The set of integers~$p$ gives the cyclic sequence of indices $\dots 123412341234\dots $ and we associate variables
$\{\tau_j\}$ as follows
\begin{gather*}
\begin{array}{@{}cccccccccccccccccc}
\dots&  1&2&3&4&1&2&3&4&1&2&3&4&1&2&3&4& \dots
\\
\dots & \tau_{-4}&\tau_{-3}&\tau_{-2}&\tau_{-1}&\tau_1&\tau_2&\tau_3&\tau_4& \tau_5&\tau_6&\tau_7& &\tau_8&\tau_9&
&\tau_{10} & \dots
\end{array}
\end{gather*}
Note that the skip between $\tau_7$ and $\tau_8$ or $\tau_9$ and $\tau_{10}$ means no corresponding variable appears in
the following crystal graph
\begin{gather}
\begin{split}
& \begin{xy}
(0,0) *{\tau_{-2}}="(-2)",
(20,0) *{\frac{\tau_{-1}\tau_2}{\tau_3}}="(-12,3)",
(40,0) *{\frac{\tau_2}{\tau_4}}="(2,4)",
(30,-10) *{\frac{\tau_{-1}\tau_5}{\tau_6}}="(-15,6)",
(50,-10) *{\frac{\tau_3\tau_5}{\tau_4\tau_6}}="(35,46)",
(70,-10) *{\frac{\tau_5}{\tau_7}}="(5,7)",
(40,-20) *{\frac{\tau_{-1}}{\tau_8}}="(-1,8)",
(60,-20) *{\frac{\tau_3}{\tau_4\tau_8}}="(3,48)",
(80,-20) *{\frac{\tau_6}{\tau_7\tau_8}}="(6,78)",
(100,-20) *{\frac{1}{\tau_9}\,\, ,}="(0,9)",
\ar@{->} "(-2)";"(-12,3)"^{3}
\ar@{->} "(-12,3)";"(2,4)"^{4}
\ar@{->} "(-12,3)";"(-15,6)"^{2}
\ar@{->} "(-15,6)";"(35,46)"^{4}
\ar@{->} "(-15,6)";"(-1,8)"^{1}
\ar@{->} "(-1,8)";"(3,48)"^{4}
\ar@{<-} "(0,9)";"(6,78)"^{2}
\ar@{<-} "(6,78)";"(5,7)"^{1}
\ar@{<-} "(5,7)";"(35,46)"_{3}
\ar@{<-} "(35,46)";"(2,4)"_{2}
\ar@{<-} "(6,78)";"(3,48)"^{3}
\ar@{<-} "(3,48)";"(35,46)"_{1}
\end{xy}
\end{split}
\label{c-gra}
\end{gather}
where the highest weight monomial $\tau_{-2}$ has a~weight $\Lambda_3$ and the lowest weight monomial $\frac{1}{\tau_9}$
has a~weight $-\Lambda_2$.

As stated above, the initial cluster variables on $G^{u,v}$ are expressed by generalized minors $\{\Delta(k;{\mathbf i})
 \,|\,  1\leq k\leq l(u)+l(v)\}$, where ${\mathbf i}$ is the reduced expression of $(u,v)\in W\times W$.
Now, as an example we consider the case $G={\rm SL}_5(\mathbb C)$ as above.
Let $W=\mathfrak{S}_5=\langle s_i \,|\, 1\leq i\leq 4\rangle$ be the symmetric group and set
$u:=s_1s_2s_3s_4s_1s_2s_3s_1s_2s_1$, $v:=e$, and set a~reduced word ${\mathbf i}$ for~$u$ as ${\mathbf
i}:=(1,2,3,4,1,2,3,1,2,1)$.
We have for $x=(x_{i,j})\in {\rm SL}_5(\mathbb C)$
\begin{gather}
\Delta(6;{\mathbf i})(x)=
\begin{vmatrix}
x_{31} & x_{32}
\\
x_{41} & x_{42}
\end{vmatrix}
=x_{31}x_{42}-x_{32}x_{41}.
\label{Del6}
\end{gather}
See Section~\ref{gmc} for the detailed explanation.

Now, let us consider the generalized minors on $L^{u,v}:=NuN \cap B_- vB_-$ instead of $G^{u,v}$ since their dif\/ference
is, indeed, only the factor from the torus part.
We call $L^{u,v}$ reduced double Bruhat cell.
For the above $u$, $v$, there exists a~birational map ${\bf x}^L_{{\mathbf i}}:(\mathbb{C}^{\times})^{10}\overset{\sim}
{\longrightarrow} L^{u,v}$ given by
\begin{gather}
{\bf x}^L_{{\bf i}}(\tau_1,\dots,\tau_{10})
\nonumber
\\
\qquad
=x_{-1}(\tau_1)x_{-2}(\tau_2)x_{-3}(\tau_3)x_{-4}(\tau_4) x_{-1}(\tau_5)x_{-2}(\tau_6)x_{-3}(\tau_7)
x_{-1}(\tau_8)x_{-2}(\tau_9) x_{-1}(\tau_{10})
\nonumber
\\
\qquad
=
\begin{pmatrix}
\dfrac{1}{\tau_1\tau_5\tau_8\tau_{10}} & 0 & 0 & 0 & 0
\\
A~& \dfrac{\tau_1\tau_5\tau_8\tau_{10}}{\tau_2\tau_6\tau_9} & 0 & 0 & 0
\\
B &C& \dfrac{\tau_2\tau_6\tau_9}{\tau_3\tau_7} & 0 & 0
\vspace{1mm}\\
D & E & \dfrac{\tau_3\tau_9}{\tau_4}+\dfrac{\tau_6\tau_9}{\tau_7} & \dfrac{\tau_3\tau_7}{\tau_4} & 0
\\
1 & \tau_{10} & \tau_9 & \tau_7 & \tau_4
\end{pmatrix},
\label{int2}
\end{gather}
\text{where}
\begin{gather*}
A=\frac{\tau_1\tau_5\tau_8}{\tau_2\tau_6\tau_9}+ \frac{\tau_1\tau_5}{\tau_2\tau_6\tau_{10}}
+\frac{\tau_1}{\tau_2\tau_8\tau_{10}} +\frac{1}{\tau_5\tau_8\tau_{10}},
\\
B=\frac{\tau_2\tau_6}{\tau_3\tau_7}+\frac{\tau_2\tau_8}{\tau_3\tau_9}
+\frac{\tau_5\tau_8}{\tau_6\tau_9}+\frac{\tau_2}{\tau_3\tau_{10}}+
\frac{\tau_5}{\tau_6\tau_{10}}+\frac{1}{\tau_8\tau_{10}},
\\
C=\frac{\tau_2\tau_6\tau_{10}}{\tau_3\tau_7} +\frac{\tau_2\tau_8\tau_{10}}{\tau_3\tau_9}
+\frac{\tau_5\tau_8\tau_{10}}{\tau_6\tau_9},
\\
D=\frac{\tau_3}{\tau_4}+\frac{\tau_6}{\tau_7} +\frac{\tau_8}{\tau_9}+\frac{1}{\tau_{10}},
\qquad
E=\frac{\tau_3\tau_{10}}{\tau_4}+\frac{\tau_6\tau_{10}}{\tau_7} +\frac{\tau_8\tau_{10}}{\tau_9},
\\
x_{-i}(t)=
\begin{matrix}
\vphantom{\ddots}\\
i^{\rm  th}\\
\\
\vphantom{\ddots}
\end{matrix}
\begin{pmatrix}
 \ddots & & & \cr
   & t^{-1} & 0 &  \cr
   & 1 & t &  \cr
   & & & \ddots
\end{pmatrix}.
\end{gather*}
Therefore, by~\eqref{Del6} and~\eqref{int2} we f\/ind
\begin{gather}
  \Delta^L(6;{\bf i})(\tau):=\big(\Delta(6;{\bf i})\circ x^L_{{\bf i}}\big)(\tau_1,\dots,\tau_{10}) =
\begin{vmatrix}
B & C
\\
D & E
\end{vmatrix}
\nonumber
\\
\hphantom{\Delta^L(6;{\bf i})(\tau)}{} =\frac{\tau_2}{\tau_4}+\frac{\tau_3\tau_5}{\tau_4\tau_6}+\frac{\tau_5}{\tau_7}
+\frac{\tau_3}{\tau_4\tau_8}+\frac{\tau_6}{\tau_7\tau_8}+\frac{1}{\tau_9}.
\label{Del-fin}
\end{gather}

Now, observing the crystal graph~\eqref{c-gra} and the Laurent polynomial~\eqref{Del-fin}, we realize that each term
in~\eqref{Del-fin} appears in~\eqref{c-gra} and they constitute so-called lower Demazure crystal associated with the
element $u_{\leq 6}\in \mathfrak{S}_5$~\cite{K3}.

Those facts motivate us to f\/ind a~new linkage between the cluster variables on $L^{u,v}\subset G^{u,v}$ and the monomial
realizations of crystals.

In this paper, we shall treat the case $G={\rm SL}_{r+1}(\mathbb{C})$, $v=e$ and some special $u\in W=\mathfrak{S}_{r+1}$.
More precisely, we treat an element $u\in W$ whose reduced word (Def\/inition~\ref{redworddef}) can be written as a~left
factor of the standard longest word $(1,2,3,\dots,r,1,2,3,\dots,(r-1),\dots,1,2,1)$:
\begin{gather*}
u=s_1s_2\cdots s_rs_1\cdots s_{r-1}\cdots s_1\cdots s_{r-m+2} s_1\cdots s_{i_n},
\end{gather*}
where $n:=l(u)$ is the length of~$u$ and $1\leq i_n\leq r-m+1$.
And we treat (reduced) double Bruhat cells of the form $G^{u,e}:=BuB \cap B_-$ and $L^{u,e}:=NuN \cap B_-$, where~$B$
(resp.\ $B_-$) is the subgroup of upper (resp.\ lower) triangular matrices in $G={\rm SL}_{r+1}(\mathbb{C})$.
Then generalized minors are a~part of classical minors (Def\/inition~\ref{Alem}).
This case matches well to the Demazure crystals.
In fact, we shall describe generalized minors in terms of summations over certain monomial realizations of Demazure
crystals in the main result Theorem~\ref{thm2}.
For example, \eqref{Del-fin}~shows that the generalized minor $\Delta^L(6;{\bf i})(\tau)$ is described in terms of
summation over certain monomial realization of the Demazure crystal $B^-_{u\leq 6}(-\Lambda_2)$.

In forthcoming paper, we shall treat more general setting, like as, the Weyl group element $v\in W$ is non-identity or
type C.~In these cases, the generalized minors are described also by monomial realizations of crystals.

\section{Factorization theorem for type A}\label{DBCs}

In this section, we shall introduce (reduced) double Bruhat cells $G^{u,v}$, $L^{u,v}$, and their properties in the case
$G={\rm SL}_{r+1}({\mathbb C})$, $v=e$ and some special $u\in W$.
In~\cite{B-Z,F-Z}, these properties had been proven for simply connected, connected, semisimple complex
algebraic groups and arbitrary $u,v\in W$.

For $l\in \mathbb{Z}_{>0}$, we set $[1,l]:=\{1,2,3,\dots,l\}$.

\subsection{Double Bruhat cells}

Let $G={\rm SL}_{r+1}(\mathbb{C})$ be the simple complex algebraic group of type $\text{A}_r$,~$B$ and $B_-$ be two opposite
Borel subgroups in~$G$, that is,~$B$ (resp.\ $B_-$) is the subgroup of upper (resp.\ lower) triangular matrices in $G={\rm SL}_{r+1}(\mathbb{C})$.
Let $N\subset B$ and $N_-\subset B_-$ be their unipotent radicals, $H:=B\cap B_-$ a~maximal torus, and $W:=\text{Norm}_G(H)/H$ the Weyl group.
In this case, Weyl group~$W$ is isomorphic to the symmetric group $\mathfrak{S}_{r+1}$.

We have two kinds of Bruhat decompositions of~$G$ as follows
\begin{gather*}
G=\displaystyle\coprod_{u \in W}BuB=\displaystyle\coprod_{u \in W}B_-uB_-.
\end{gather*}
Then, for~$u$, $v\in W$, we def\/ine the {\it double Bruhat cell} $G^{u,v}$ as follows
\begin{gather*}
G^{u,v}:=BuB \cap B_-vB_-.
\end{gather*}
This is biregularly isomorphic to a~Zariski open subset of an af\/f\/ine space of dimension $r+l(u)+l(v)$~\cite[Theorem~1.1]{F-Z}.

We also def\/ine the {\it reduced double Bruhat cell} $L^{u,v}$ as follows
\begin{gather*}
L^{u,v}:=NuN \cap B_-vB_- \subset G^{u,v}.
\end{gather*}
As is similar to the case $G^{u,v}$, $L^{u,v}$ is biregularly isomorphic to a~Zariski open subset of an af\/f\/ine space of
dimension $l(u)+l(v)$~\cite[Proposition 4.4]{B-Z}.

\begin{defn}
\label{redworddef}
Let $u=s_{i_1}\cdots s_{i_n}$ be a~reduced expression of $u\in W$, $i_1,\dots,i_n\in [1,r]$.
Then the f\/inite sequence
\begin{gather*}
{\bf i}:=(i_1,\dots,i_n)
\end{gather*}
is called {\it reduced word} ${\bf i}$ for~$u$.
\end{defn}

In this paper, we treat (reduced) double Bruhat cells of the form $G^{u,e}:=BuB \cap B_-$ and $L^{u,e}:=NuN \cap B_-$,
where $u\in W$ is an element whose reduced word can be written as a~left factor of
$(1,2,3,\dots,r,1,2,3,\dots,(r-1),\dots,1,2,1)$:
\begin{gather}
\label{uset0}
u=s_1s_2\cdots s_rs_1\cdots s_{r-1}\cdots s_1\cdots s_{r-m+2} s_1\cdots s_{i_n},
\end{gather}
where $n:=l(u)$ is the length of~$u$ and $1\leq i_n\leq r-m+1$.
Let ${\bf i}$ be a~reduced word of~$u$:
\begin{gather}
\label{iset0}
{\bf i}=(\underbrace{1,\dots,r}_{\text{1$^{\rm st}$ cycle}},\underbrace{1,\dots,(r-1)}_{\text{2$^{\rm nd}$
cycle}},\dots,\underbrace{1,\dots,(r-m+2)}_{\text{$(m-1)^{\rm th}$ cycle}}, \underbrace{1,\dots,i_n}_{\text{$m^{\rm th}$ cycle}}).
\end{gather}
Note that $(1,2,3,\dots,r,1,2,3,\dots,(r-1),\dots,1,2,1)$ is a~reduced word of the longest element in~$W$.

\subsection{Factorization theorem for type A}%\label{factproA}

In this subsection, we shall introduce the isomorphisms between double Bruhat cell $G^{u,e}$ and $H\times
(\mathbb{C}^{\times})^{l(u)}$, and between $L^{u,e}$ and $(\mathbb{C}^{\times})^{l(u)}$.
As in the previous section, we consider the case $G:={\rm SL}_{r+1}(\mathbb{C})$.
We set $\mathfrak g:=\text{Lie}(G)$ with the Cartan decomposition $\mathfrak g={\mathfrak n}_-\oplus {\mathfrak h}\oplus{\mathfrak n}$.
Let~$e_i$,~$f_i$ ($i\in[1,r]$)
be the generators of ${\mathfrak n}$, ${\mathfrak n}_-$.
For $i\in[1,r]$ and $t \in \mathbb{C}$, we set $x_i(t):=\exp(te_i)$, $y_{i}:=\exp(tf_i)$.
Let $\varphi_i:{\rm SL}_2(\mathbb{C}) \rightarrow G$ be the canonical embedding corresponding to each simple root $\alpha_i$.
Then we have
\begin{gather*}
x_i(t)=\varphi_i
\begin{pmatrix}
1 & t
\\
0 & 1
\end{pmatrix},
\qquad
y_{i}(t)=\varphi_i
\begin{pmatrix}
1 & 0
\\
t & 1
\end{pmatrix}.
\end{gather*}
We can express $x_i(t)$, $y_{i}(t)$, as the following matrices
\begin{gather}
x_i(t)=
\begin{matrix}
\vphantom{\ddots}\\
i^{\rm  th}\\
\\
\vphantom{\ddots}
\end{matrix}
\begin{pmatrix}
  \ddots & & &  \cr
  & 1 & t & \cr
  & 0 & 1 & \cr
  & & & \ddots
\end{pmatrix},
\qquad
y_{i}(t)=
\begin{matrix}
\vphantom{\ddots}\\
i^{\rm  th}\\
\\
\vphantom{\ddots}
\end{matrix}
\begin{pmatrix}
   \ddots & & &  \cr
  & 1 & 0 & \cr
  & t & 1 & \cr
  & & & \ddots
\end{pmatrix}.
\label{transmat0}
\end{gather}
For a~reduced word ${\bf i}=(i_1,i_2,\dots,i_n)$, we def\/ine a~map $x^G_{{\bf i}}:H\times \mathbb{C}^n \rightarrow
G$ as
\begin{gather*}
x^G_{{\bf i}}(a; t_1, \dots t_n):=a\cdot y_{i_1}(t_1)y_{i_2}(t_2)\cdots y_{i_n}(t_n).
\end{gather*}

\begin{thm}[\protect{\cite[Theorem~1.2]{F-Z}}]\label{fp}
We set $u\in W$ and its reduced word ${\bf i}$ as in~\eqref{uset0} and \eqref{iset0}.
The map $x^G_{{\bf i}}$ defined above can be restricted to a~biregular isomorphism between $H\times
(\mathbb{C}^{\times})^{l(u)}$ and a~Zariski open subset of $G^{u,e}$.
\end{thm}

Next, for $i \in [1,r]$ and $t\in \mathbb{C}^{\times}$, we def\/ine as follows
\begin{gather*}
\alpha_i^{\vee}(t):=\varphi_i
\begin{pmatrix}
t & 0
\\
0 & t^{-1}
\end{pmatrix},
\qquad x_{-i}(t):=y_{i}(t)\alpha_i^{\vee}(t^{-1})=\varphi_i
\begin{pmatrix}
t^{-1} & 0
\\
1 & t
\end{pmatrix}.
\end{gather*}
We can express $x_{-i}(t)$ and $\alpha_i^{\vee}(t)$ as the following matrices
\begin{gather}
\label{transmat}
x_{-i}(t)=
\begin{matrix}
\vphantom{\ddots}\\
i^{\rm  th}\\
\\
\vphantom{\ddots}
\end{matrix}
\begin{pmatrix}
 \ddots & & &  \cr
 & t^{-1} & 0 & \cr
 & 1 & t & \cr
 & & & \ddots
\end{pmatrix},
\qquad \alpha_i^{\vee}(t)=\diag(1,\dots,1,\stackrel{i}{\stackrel{\vee}{t}}, \stackrel{i+1}{\stackrel{\vee}{t^{-1}}},1,\dots,1).
\end{gather}

For ${\bf i}=(i_1, \dots,i_n)$ ($i_1,\dots,i_n\in[1,r]$), we def\/ine a~map $x^L_{{\bf i}}:\mathbb{C}^n \rightarrow
G$ as
\begin{gather*}
x^L_{{\bf i}}(t_1, \dots, t_n):=x_{-i_1}(t_1)\cdots x_{-i_n}(t_n).
\end{gather*}
We have the following theorem which is similar to the previous one.

\begin{thm}[\protect{\cite[Proposition~4.5]{B-Z}}]  \label{fp2}
We set $u\in W$ and its reduced word ${\bf i}$ as in~\eqref{uset0} and~\eqref{iset0}.
The map $x^L_{{\bf i}}$ defined above can be restricted to a~biregular isomorphism between $
(\mathbb{C}^{\times})^{l(u)}$ and a~Zariski open subset of $L^{u,e}$.
\end{thm}

Finally, we def\/ine a~map $\bar{x}^G_{{\bf i}}:H\times(\mathbb{C}^{\times})^{n}\rightarrow G^{u,v}$ as
\begin{gather*}
\bar{x}^G_{{\bf i}}(a;t_1,\dots,t_n) =ax^L_{{\bf i}}(t_1,\dots,t_n).
\end{gather*}
\begin{prop}
\label{gprime}
In the above setting, the map $\bar{x}^G_{{\bf i}}$ is a~biregular isomorphism between
$H\times(\mathbb{C}^{\times})^{n}$ and a~Zariski open subset of $G^{u,e}$.
\end{prop}

\begin{proof}
We set $l_0:=0$, $l_1:=r$, $l_2:=r+(r-1), \dots, l_m:=r+(r-1)+\dots+(r-m+1)$.
We def\/ine a~map $\phi:H\times(\mathbb{C}^{\times})^{n}\rightarrow H\times(\mathbb{C}^{\times})^{n}$,
${\bf t}=(a;t_1,\dots,t_n) \mapsto(a({\bf t});\tau_1({\bf t}),\dots,\tau_n({\bf t}))$ as
\begin{gather}
a({\bf t})=a\cdot\underbrace{\alpha_1^{\vee}(t_1)^{-1}\cdots \alpha_r^{\vee}(t_r)^{-1}}_{\text{$1^{\rm st}$ cycle}}\cdots
\underbrace{\alpha_1^{\vee}(t_{l_{m-1}+1})^{-1}\cdots \alpha_{i_n}^{\vee}(t_{l_{m-1}+i_n})^{-1}}_{\text{$m^{\rm th}$ cycle}},
\nonumber
\\
\label{mbase0}
\tau_{l_{s}+j}({\bf t})=\frac{(t_{l_{s+1}+j-1}t_{l_{s+2}+j-1}\cdots t_{l_{m-1}+j-1})(t_{l_{s}+j+1}t_{l_{s+1}+j+1}\cdots
t_{l_{m-1}+j+1})}{t_{l_{s}+j}(t_{l_{s+1}+j}\cdots t_{l_{m-1}+j})^{2}},
\end{gather}
where in~\eqref{mbase0}, if {\bf i} does not include~$j$ (resp.~$j+1$, $j-1$) in~$\zeta^{\rm th}$ cycle then we set
$t_{l_{\zeta-1}+j}=1$ (resp.~$t_{l_{\zeta-1}+j+1}=1$, $t_{l_{\zeta-1}+j-1}=1$).
This is a~biregular isomorphism.

Let us prove
\begin{gather*}
\bar{x}^G_{{\bf i}}(a;t_1,\dots,t_n)=\big(x^G_{{\bf i}}\circ\phi\big)(a;t_1,\dots,t_n),
\end{gather*}
which implies that $\bar{x}^G_{{\bf i}}:H\times(\mathbb{C}^{\times})^{n}\rightarrow G^{u,e}$ is a~biregular isomorphism
by Theorem~\ref{fp}.

First, we can verify the following relations by the explicit forms~\eqref{transmat0},~\eqref{transmat} and direct
calculations:
\begin{gather}
\label{base2}
\alpha_i^{\vee}(c)^{-1}y_{j}(t)=
\begin{cases}
y_{i}(c^2t)\alpha_i^{\vee}(c)^{-1} & \text{if} \quad  i=j,
\\
y_{j}(c^{-1}t)\alpha_i^{\vee}(c)^{-1} & \text{if} \quad  |i-j|=1,
\\
y_{j}(t)\alpha_i^{\vee}(c)^{-1} & \text{otherwise},
\end{cases}
\end{gather}
for $1\leq i$, $j\leq r$ and $c,t\in \mathbb{C}^{\times}$.

On the other hand, we obtain
\begin{gather}
(x^G_{{\bf i}}\circ\phi)(a;t_1,\dots,t_n) =a\times\alpha_1^{\vee}(t_1)^{-1}\cdots \alpha_r^{\vee}(t_r)^{-1}\cdots
\alpha_1^{\vee}(t_{l_{m-1}+1})^{-1}\cdots\alpha_{i_n}^{\vee}(t_{l_{m-1}+i_n})^{-1}
\nonumber
\\
\qquad{}
\times
y_{1}(\tau_1({\bf t}))y_{2}(\tau_2({\bf t}))\cdots y_{r}(\tau_r({\bf t}))\cdots
y_{1}(\tau_{l_{m-1}+1}({\bf t})) \cdots y_{i_n}(\tau_{l_{m-1}+i_n}({\bf t})).
\label{xigp}
\end{gather}

For each~$s$ and~$j$, let us move $\alpha_j^{\vee}(t_{l_s+j})^{-1}$,
$\alpha_{j+1}^{\vee}(t_{l_s+j+1})^{-1},\dots, \alpha_{i_n}^{\vee}(t_{l_{m-1}+i_n})^{-1}$
to the right of $y_{j}(\tau_{l_s+j}({\bf t}))$ by using the relations~\eqref{base2}.
For example,
\begin{gather*}
\alpha_j^{\vee}(t_{l_s+j})^{-1}\cdots\alpha_{j-1}^{\vee}(t_{l_{m-1}+j-1})^{-1}\alpha_{j}^{\vee}(t_{l_{m-1}+j})^{-1}
\\
\qquad
\phantom{=}{}
\times
\alpha_{j+1}^{\vee}(t_{l_{m-1}+j+1})^{-1}\cdots \alpha_{i_n}^{\vee}(t_{l_{m-1}+i_n})^{-1}y_{j}(\tau_{l_s+j}({\bf t}))
\\
\qquad
=\alpha_j^{\vee}(t_{l_s+j})^{-1}\cdots y_{j}\left(\frac{t_{l_{m-1}+j}^2}{t_{l_{m-1}+j-1}t_{l_{m-1}+j+1}}
\tau_{l_s+j}({\bf t})\right) \alpha_{j-1}^{\vee}(t_{l_{m-1}+j-1})^{-1}
\\
\qquad
\phantom{=}{}
\times
\alpha_{j}^{\vee}(t_{l_{m-1}+j})^{-1}\alpha_{j+1}^{\vee}(t_{l_{m-1}+j+1})^{-1}\cdots
\alpha_{i_n}^{\vee}(t_{l_{m-1}+i_n})^{-1}.
\end{gather*}

Repeating this argument, we have
\begin{gather*}
=
y_{j}\left(\!\frac{(t_{l_{s}+j}\cdots t_{l_{m-1}+j})^2}{(t_{l_{s}+j-1}\cdots t_{l_{m-1}+j-1})(t_{l_{s}+j+1}\cdots
t_{l_{m-1}+j+1})}\tau_{l_s+j}({\bf t})\!\right) \alpha_j^{\vee}(t_{l_s+j})^{-1}\cdots
\alpha_{i_n}^{\vee}(t_{l_{m-1}+i_n})^{-1}.
\end{gather*}
Note that $\frac{(t_{l_{s}+j}\dots t_{l_{m-1}+j})^2}{(t_{l_{s}+j-1}\cdots t_{l_{m-1}+j-1})(t_{l_{s}+j+1}\cdots
t_{l_{m-1}+j+1})} \tau_{l_s+j}({\bf t})=t_{l_{s}+j}$.
By~\eqref{xigp}, we have
\begin{gather*}
  \big(x^G_{{\bf i}}\circ\phi\big)(a;t_1,\dots,t_n)=a\cdot y_{1}(t_1)\alpha_1^{\vee}(t_1)^{-1}\cdots
y_{r}(t_r)\alpha_r^{\vee}(t_r)^{-1}\cdots
\\
\qquad
\phantom{=}{}
\times
y_{1}(t_{l_{m-1}+1})\alpha_1^{\vee}(t_{l_{m-1}+1})^{-1}\cdots
y_{i_n}(t_{l_{m-1}+i_n})\alpha_{i_n}^{\vee}(t_{l_{m-1}+i_n})^{-1}
\\
\qquad{}
 =a\cdot x_{-1}(t_1)\cdots x_{-r}(t_r)\cdots x_{-1}(t_{l_{m-1}+1})\cdots x_{-i_n}(t_{l_{m-1}+i_n})
=\bar{x}^G_{{\bf i}}(a;t_1,\dots,t_n).
\tag*{\qed} \end{gather*}\renewcommand{\qed}{} \end{proof}

\section{Monomial realizations of crystal bases}\label{section3}

In this section, we review the monomial realizations of crystals~\cite{K2,K, Nj}.
Let $I:=\{1,2,\dots,r\}$ be a~f\/inite index set.

\subsection{Monomial realizations of crystal bases for type A}\label{monoreal}

\begin{defn}
Let $A=(a_{ij})_{i,j\in I}$ be the Cartan matrix of type A$_r$: $A=(a_{ij})_{i,j\in I}$ is def\/ined~as
\begin{gather}
\label{carmatA}
a_{ij}=
\begin{cases}
2 & \text{if}\quad  i=j,
\\
-1 & \text{if}\quad  |i-j|=1,
\\
0 & \text{otherwise}.
\end{cases}
\end{gather}
Let $\Pi=\{\alpha_i \,|\, i\in I\}$ (resp.\ $\Pi^{\vee}=\{h_i \,|\, i\in I\}$) be the set of simple roots (resp.\ co-roots),
and~$P$ be the weight lattice.
A~{\it crystal} associated with the Cartan matrix~$A$ is a~set~$B$ together with the maps $\text{wt}: B \rightarrow P$,
$\tilde{e_{i}}$, $\tilde{f_{i}}: B \cup \{0\} \rightarrow B \cup \{0\}$ and $\varepsilon_i$,
$\varphi_i: B \rightarrow {\mathbb Z} \cup \{-\infty \}$, $i \in I$, satisfying the following properties: For $b\in B$, $i \in I$,
\begin{enumerate}\itemsep=0pt
\item[(i)] $\varphi_i(b) - \varepsilon_i(b) = \langle h_i$, $\text{wt}(b)\rangle$,
\item[(ii)] $\text{wt}(\tilde{e_{i}}b)=\text{wt}(b) + \alpha_i$, if $\tilde{e_{i}}b \in B$,
\item[(iii)] $\text{wt}(\tilde{f_{i}}b)=\text{wt}(b) - \alpha_i$, if $\tilde{f_{i}}b \in B$,
\item[(iv)] $\varepsilon_i(\tilde{e_{i}}b)=\varepsilon_i(b) -1$, $\varphi_i(\tilde{e_{i}}b)=\varphi_i(b) +1$ if
$\tilde{e_{i}}b \in B$,
\item[(v)] $\varepsilon_i(\tilde{f_{i}}b)=\varepsilon_i(b) +1$, $\varphi_i(\tilde{f_{i}}b)=\varphi_i(b) -1$ if
$\tilde{f_{i}}b \in B$,
\item[(vi)] $\tilde{f_{i}}b=b' \Leftrightarrow b=\tilde{e_{i}}b'$, if $b$, $b' \in B$,
\item[(vii)] $\varphi_i(b)=-\infty$, $b \in B$, $\Rightarrow \tilde{e_{i}}b=\tilde{f_{i}}b=0$.
\end{enumerate}
\end{defn}

Let $U_q(\mathfrak g)$ be the universal enveloping algebra associated with the Cartan matrix~$A$ in~\eqref{carmatA}, and
$\mathfrak g={\mathfrak s}{\mathfrak l}_{r+1}(\mathbb{C})$.
Let $B^{+}(\lambda)$ (resp.\ $B^{-}(\lambda)$) be the crystal base of the $U_q(\mathfrak g)$-highest (resp.\ lowest) weight module~\cite{H-K,KN}.
Note that $B^+(\lambda)=B^-(w_0\lambda)$, where $w_0$ is the longest element of~$W$.
In particular, in the case $\lambda=M\Lambda_d$, $M\in \mathbb{Z}_{>0}$, we have
\begin{gather*}
B^{+}(M\Lambda_d)=B^-(-M\Lambda_{r-d+1}).
\end{gather*}

Let us introduce monomial realizations which realize each element of $B^{\pm}(\lambda)$ as a~certain Laurent monomial.

First, we def\/ine a~set of integers $p=(p_{j,i})_{j,i \in I,\; j \neq i}$ such that
\begin{gather*}
p_{j,i}=
\begin{cases}
1 & \text{if} \quad  j<i,
\\
0 & \text{if} \quad  i<j.
\end{cases}
\end{gather*}
Second, for doubly-indexed variables $\{Y_{s,i} \,|\, i \in I$, $s\in \mathbb{Z}\}$, we def\/ine the set of monomials
\begin{gather}
\label{cydef}
{\mathcal Y}:=\left\{Y=\prod\limits_{s \in \mathbb{Z},\ i \in I}
Y_{s,i}^{\zeta_{s,i}}\, \Bigg| \,\zeta_{s,i} \in \mathbb{Z},\
\zeta_{s,i} =0~\text{except for f\/initely many}~(s,i) \right\}.
\end{gather}

Finally, we def\/ine maps $\text{wt}: {\mathcal Y} \rightarrow P$, $\varepsilon_i$, $\varphi_i: {\mathcal Y} \rightarrow
\mathbb{Z}$, $i \in I$. %, as follows.
For $Y=\prod\limits_{s \in \mathbb{Z},\; i \in I} Y_{s,i}^{\zeta_{s,i}}\in {\mathcal Y}$,
\begin{gather*}
\text{wt}(Y):= \sum\limits_{i,s}\zeta_{s,i}\Lambda_i,\!
\qquad
\varphi_i(Y):=\max\left\{\! \sum\limits_{k\leq s}\zeta_{k,i}  \,|\, s\in \mathbb{Z} \!\right\},\!
\qquad
\varepsilon_i(Y):=\varphi_i(Y)-\text{wt}(Y)(h_i).
\end{gather*}

We set
\begin{gather}
\label{asidef}
A_{s,i}:=Y_{s,i}Y_{s+1,i}\prod\limits_{j\neq i}Y_{s+p_{j,i},j}^{a_{j,i}}=
\begin{cases}
\dfrac{Y_{s,1}Y_{s+1,1}}{Y_{s,2}} & \text{if}\quad  i=1,
\vspace{1mm}\\
\dfrac{Y_{s,i}Y_{s+1,i}}{Y_{s,i+1}Y_{s+1,i-1}} & \text{if}\quad  2\leq i\leq r-1,
\vspace{1mm}\\
\dfrac{Y_{s,r}Y_{s+1,r}}{Y_{s+1,r-1}} & \text{if}\quad  i=r,
\end{cases}
\end{gather}
and def\/ine the Kashiwara operators as follows
\begin{gather*}
\tilde{f}_iY=
\begin{cases}
A_{n_{f_i},i}^{-1}Y & \text{if} \quad  \varphi_i(Y)>0,
\\
0 & \text{if} \quad  \varphi_i(Y)=0,
\end{cases}
\qquad
\tilde{e}_iY=
\begin{cases}
A_{n_{e_i},i}Y & \text{if} \quad  \varepsilon_i(Y)>0,
\\
0 & \text{if} \quad  \varepsilon_i(Y)=0,
\end{cases}
\end{gather*}
where
\begin{gather}
\label{nxi}
n_{f_i}:=\min \left\{n \,\Bigg|\, \varphi_i(Y)= \sum\limits_{k\leq n}\zeta_{k,i}\right\},
\qquad
n_{e_i}:=\max \left\{n \,\Bigg|\, \varphi_i(Y)= \sum\limits_{k\leq n}\zeta_{k,i}\right\}.
\end{gather}

Then the following theorem holds:

\begin{thm}[\cite{K, Nj}]\label{monorealmain}\quad
\begin{enumerate}\itemsep=0pt
\item[(i)] For the set $p=(p_{j,i})$ as above, $({\mathcal Y}, \text{wt}, \varphi_i, \varepsilon_i,\tilde{f}_i,
\tilde{e}_i)_{i\in I}$ is a~crystal.
When we emphasize~$p$, we write ${\mathcal Y}$ as ${\mathcal Y}(p)$.
\item[(ii)] If a~monomial $Y \in {\mathcal Y}(p)$ satisfies $\varepsilon_i(Y)=0$ $($resp.~$\varphi_i(Y)=0)$ for all $i \in I$,
then the connected component containing~$Y$ is isomorphic to $B^+(\text{wt}(Y))$ $($resp.\ $B^-(\text{wt}(Y)))$.
\end{enumerate}
\end{thm}

\begin{defn}
\label{yemb}
Let $Y\in{\mathcal Y}(p)$ be a~monomial and let $\mathbb B$ be the unique connected component in ${\mathcal Y}(p)$
including~$Y$.
Suppose that~$\lambda$ is the highest (resp.\ lowest) weight of $\mathbb B$.
We denote the embedding
\begin{gather*}
\mu_{Y}: \ B^+(\lambda)\hookrightarrow\mathbb B\subset {\mathcal Y}(p)
\qquad
(\text{resp.} \
\mu_{Y}: \ B^-(\lambda)\hookrightarrow\mathbb B).
\end{gather*}
\end{defn}

Note that if~$Y$ and $Y'$ are in the same component then $\mu_Y=\mu_{Y'}$.

\begin{rem}
\label{actionrem}
The actions of $\tilde{e}_i$ and $\tilde{f}_i$ on~$Y$ are determined by wt$(Y)(h_i)$, $\varphi_i(Y)$ and
$\varepsilon_i(Y)$, which are determined by the factors~$Y^{\pm1}_{s,i}$, $s\in \mathbb{Z}$.
Thus, when we consider the actions of~$\tilde{e}_i$ and~$\tilde{f}_i$, we need to see the factors
$\{Y^{\pm1}_{s,i}\}_{s\in\mathbb{Z}}$ only.
\end{rem}

\begin{ex}
\label{emb}
For $\lambda=\beta \Lambda_d$ (resp.~$\lambda=-\beta \Lambda_d$), $\beta\in \mathbb{Z}_{>0}$, $d\in I$,
we can embed $B^{+}(\lambda)$ (resp.~$B^{-}(\lambda)$) in ${\mathcal Y}$ as a~crystal~by
\begin{gather*}
v_{\lambda} \mapsto Y_{\beta+\gamma,d}Y_{\beta-1+\gamma,d}\cdots Y_{1+\gamma,d},
\qquad
\left(\text{resp.} \ v_{\lambda} \mapsto
\frac{1}{Y_{\beta+\gamma,d}Y_{\beta-1+\gamma,d}\cdots Y_{1+\gamma,d}}\right),
\end{gather*}
where $v_{\lambda}$ is the highest (resp.\ lowest) weight vector of $B^{+}(\lambda)$ (resp.~$B^{-}(\lambda)$),
and~$\gamma$ is an arbitrary integer.
For $Y^+:=Y_{\beta+\gamma,d}Y_{\beta-1+\gamma,d}\cdots Y_{1+\gamma,d}$
(resp.~$Y^-:=\frac{1}{Y_{\beta+\gamma,d}Y_{\beta-1+\gamma,d}\cdots Y_{1+\gamma,d}}$), $\mu_{Y^+}$ (resp.~$\mu_{Y^-}$) denotes
the embedding in Definition~\ref{yemb}.
Then $Y^+$ (resp.~$Y^-$) is the highest (resp.\ lowest) weight vector in $\mu_{Y^+}(B^{+}(\lambda))$ (resp.~$\mu_{Y^-}(B^{-}(\lambda))$).
\end{ex}

We set $l_0:=0$, $l_1:=r$, $l_2:=r+(r-1), \dots, l_s:=r+(r-1)+\cdots+(r-s+1), \dots$, $l_r:=r+(r-1)+\cdots+2+1$ and changing
the variables $Y_{s,j}$ to $\tau_{l_{s}+j}$, $1\leq j\leq r-s$.
For $s<0$, we transform the variables $Y_{s,j}$ to $\tau_{-(r+1-j)}$, $1\leq j\leq r$,
\begin{gather*}
\begin{array}
{cccccccccccccccc} \dots&r&1&\dots&r-1&r&1&\dots&r-1&r&1&2&\dots
\\
\dots&\tau_{-1} &\tau_{l_0+1}&\dots&\tau_{l_0+r-1}&\tau_{l_0+r}& \tau_{l_1+1}&\dots&\tau_{l_1+r-1}&
&\tau_{l_2+1}&\tau_{l_2+2}&\dots
\end{array}
\end{gather*}

\begin{rem}
In the above setting, the variables $\{Y_{s,j} \,|\, r-s<j\}$ do not correspond to any variables in~$\tau$.
As we have seen in~\eqref{c-gra}, these variables do not appear in the crystal base which we treat in this paper.
In other words, we only need variables associated with
\begin{gather*}
{\bf j}=(1,\dots,r,1,\dots,r-1,\dots,1,2,1),
\end{gather*}
which coincides with a~specific reduced word of the longest element of~$W$.
\end{rem}
\begin{rem}
\label{tautau}
For the variables $\tau_{l_s+{0}}$, $\tau_{l_s+{r+1}}$ $ (0\leq s\leq m-1)$
we understand
\begin{gather*}
\tau_{l_s+{0}}=\tau_{l_s+{r+1}}=1.
\end{gather*}
For example, if $i=1$ then
\begin{gather*}
\tau_{l_s+i-1}=1.
\end{gather*}
\end{rem}
\begin{ex}
Let us consider the action of $\tilde{e}_1$ on the monomial $\frac{1}{\tau_{l_{r-1}+1}}$.
Following the method in Section~\ref{monoreal}, we have $\text{wt}(\frac{1}{\tau_{l_{r-1}+1}})=-\Lambda_1$,
$\varphi_1(\frac{1}{\tau_{l_{r-1}+1}})=0$, $\varepsilon_1(\frac{1}{\tau_{l_{r-1}+1}})
=\varphi_1(\frac{1}{\tau_{l_{r-1}+1}})-\text{wt}(\frac{1}{\tau_{l_{r-1}+1}})(h_1)=1$, and $n_{e_1}=r-2$.
Thus, since we have $A_{r-2,1}=\tau_{l_{r-2}+1}\tau_{l_{r-1}+1}\tau_{l_{r-2+p_{2,1}}+2}^{a_{2,1}}
=\frac{\tau_{l_{r-2}+1}\tau_{l_{r-1}+1}}{\tau_{l_{r-2}+2}}$, we get
\begin{gather*}
\tilde{e}_1\frac{1}{\tau_{l_{r-1}+1}}=A_{r-2,1}\frac{1}{\tau_{l_{r-1}+1}} =\frac{\tau_{l_{r-2}+1}}{\tau_{l_{r-2}+2}}.
\end{gather*}
Similarly, we have
\begin{gather*}
\tilde{e}_2\tilde{e}_1\frac{1}{\tau_{l_{r-1}+1}} =A_{r-3,2}\frac{\tau_{l_{r-2}+1}}{\tau_{l_{r-2}+2}}=\frac{\tau_{l_{r-3}+2}}{\tau_{l_{r-3}+3}},
\qquad
 A_{r-3,2}=\frac{\tau_{l_{r-3}+2}\tau_{l_{r-2}+2}}{\tau_{l_{r-3}+3}\tau_{l_{r-2}+1}} ,
\\
\tilde{e}_3\tilde{e}_2\tilde{e}_1\frac{1}{\tau_{l_{r-1}+1}}
=A_{r-4,3}\frac{\tau_{l_{r-3}+2}}{\tau_{l_{r-3}+3}}
=\frac{\tau_{l_{r-4}+3}}{\tau_{l_{r-4}+4}},
\qquad
 A_{r-4,3}=\frac{\tau_{l_{r-4}+3}\tau_{l_{r-3}+3}}{\tau_{l_{r-4}+4}\tau_{l_{r-3}+2}} .
\end{gather*}
Applying $\tilde{e}_i$ repeatedly, we obtain
\begin{gather*}
\tilde{e}_k\cdots\tilde{e}_2\tilde{e}_1\frac{1}{\tau_{l_{r-1}+1}}
=A_{r-1-k,k}\tilde{e}_{k-1}\cdots\tilde{e}_2\tilde{e}_1\frac{1}{\tau_{l_{r-1}+1}}
=\frac{\tau_{l_{r-1-k}+k}}{\tau_{l_{r-1-k}+k+1}},
\qquad
k=1,\dots,r,
\end{gather*}
where, $A_{r-1-k,k}=\frac{\tau_{l_{r-1-k}+k} \tau_{l_{r-k}+k}}{\tau_{l_{r-1-k}+k+1}\tau_{l_{r-k}+k-1}}$,
$\tau_{l_{-1}+r}=r$, $\tau_{l_{-1}+r+1}:=1$.
For $i\in I$, we have $\varphi_i(\frac{1}{\tau_{l_{r-1}+1}})=0$.
Hence $\tilde{f}_i(\frac{1}{\tau_{l_{r-1}+1}})=0$.
\end{ex}
\begin{ex}
\label{exex}
For a~given $i\in I$ and $Y=\prod\limits_{s \in \mathbb{Z}} \tau_{l_s+i}^{\zeta_{s,i}}$, we def\/ine
$\nu_Y(n):=\sum\limits_{s\leq n}\zeta_{s,i}$.

For $j\in \mathbb{Z}_{>0}$, we set
\begin{gather*}
Y=\frac{1}{\tau_{l_{q_1}+i}\tau_{l_{q_2}+i}\cdots\tau_{l_{q_j}+i}},
\qquad
0\leq q_1<q_2<\dots<q_j\leq r-1.
\end{gather*}

First, let us calculate $n_{e_i}$ \eqref{nxi}.
We obtain wt$(Y)=-j\Lambda_i$ and
\begin{gather*}
\nu_Y(n)=0
\qquad
\text{for}
\quad
n<0,
\\
\nu_Y(0)=\nu_Y(1)=\cdots=\nu_Y(q_1-1)=0,
\qquad
\nu_Y(q_1)=\nu_Y(q_1+1)=\cdots=\nu_Y(q_2-1)=-1,
\\
\nu_Y(q_2)=\nu_Y(q_2+1)=\cdots=\nu_Y(q_3-1)=-2,
\qquad
\nu_Y(q_3)=\dots=\nu_Y(q_4-1)=-3,  \ \ \dots.
\end{gather*}
Thus, we get $\varphi_i(Y)=\max \{\nu_Y(n) \,|\, n\in \mathbb{Z} \}=0$ and
\begin{gather*}
n_{e_i}=\max \{n \,|\, \nu_Y(n)=0\}=q_1-1.
\end{gather*}
Next, since $\text{wt}(Y)(h_i)=-j$, we have $\varepsilon_i(Y)=\varphi_i(Y)-\text{wt}(Y)(h_i)=j>0$.
Therefore,
\begin{gather*}
\tilde{e}_i Y=A_{q_1-1,i}Y=\frac{\tau_{l_{q_1-1}+i}}{\tau_{l_{q_1-1}+i+1}\tau_{l_{q_1}+i-1}\tau_{l_{q_2}+i}\cdots\tau_{l_{q_j}+i}},
\qquad A_{q_1-1,i}=\frac{\tau_{l_{q_1-1}+i}\tau_{l_{q_1}+i}}{\tau_{l_{q_1-1}+i+1}\tau_{l_{q_1}+i-1}}  .
\end{gather*}
Similarly, for $k=1,2,\dots,j$, we get
\begin{gather*}
\tilde{e}^k_i Y=A_{q_k-1,i}\cdots A_{q_2-1,i}A_{q_1-1,i}Y=\prod\limits^k_{s=1}
\left(\frac{\tau_{l_{q_k-1}+i}}{\tau_{l_{q_k-1}+i+1}\tau_{l_{q_k}+i-1}}\right)
\frac{1}{\tau_{l_{q_{k+1}}+i}\cdots\tau_{l_{q_j}+i}}.
\end{gather*}
\end{ex}

\subsection{Demazure crystal}

For $w\in W$, let us def\/ine an {\it upper Demazure crystal} $B^{+}_{w}(\lambda)$.
This is a~subset of the crystal $B^{+}(\lambda)$ def\/ined as follows.

\begin{defn}
Let $u_{\lambda}$ be the highest weight vector of $B^{+}(\lambda)$.
For the identity element~$e$ of~$W$, we set $B^{+}_{e}(\lambda):=\{u_{\lambda}\}$.
For $w\in W$, if $s_iw<w$,
\begin{gather*}
B^{+}_{w}(\lambda):=\big\{\tilde{f}_i^kb \,|\, k\geq0,\; b\in B^{+}_{s_iw}(\lambda),\; \tilde{e}_ib=0\big\}\setminus \{0\}.
\end{gather*}
\end{defn}
Similarly, we def\/ine a~{\it lower Demazure crystal} $B^{-}_{w}(\lambda)$ as follows.

\begin{defn}
Let $v_{\lambda}$ be the lowest weight vector of $B^{-}(\lambda)$.
We set $B^{-}_{e}(\lambda):=\{v_{\lambda}\}$.
For $w\in W$, if $s_iw<w$,
\begin{gather*}
B^{-}_{w}(\lambda):=\big\{\tilde{e}_i^kb \,|\, k\geq0,\; b\in B^{-}_{s_iw}(\lambda),\; \tilde{f}_ib=0\big\}\setminus \{0\}.
\end{gather*}
\end{defn}
\begin{thm}[\cite{K3}]\label{kashidem}
For $w\in W$, let $w=s_{i_1}\cdots s_{i_n}$ be an arbitrary reduced expression.
Let $u_{\lambda}$ $($resp.~$v_{\lambda'})$ be the highest $($resp.\ lowest$)$ weight vector of $B^{+}(\lambda)$ $($resp.~$B^{-}(\lambda'))$.
Then
\begin{gather*}
  B^{+}_{w}(\lambda)=\big\{\tilde{f}_{i_1}^{a(1)}\cdots
\tilde{f}_{i_n}^{a(n)}u_{\lambda} \,|\, a(1),\dots,a(n)\in\mathbb{Z}_{\geq0} \big\} \setminus\{0\},
\\
 B^{-}_{w}(\lambda')=\big\{\tilde{e}_{i_1}^{a(1)}\cdots
\tilde{e}_{i_n}^{a(n)}v_{\lambda'} \,|\, a(1),\dots,a(n)\in\mathbb{Z}_{\geq0} \big\}\setminus\{0\}.
\end{gather*}
\end{thm}
Let $P^+$ be the set of dominant weights.
We set $P^-:=-P^+$.
\begin{defn}
\label{dempoly}
Let ${\mathcal Y}(p)$ be the monomial realization of crystal associated with $p=(p_{i,j})$.
Suppose that $Y\in {\mathcal Y}(p)$ be a~highest (resp.\ lowest) monomial with a~weight $\lambda\in P^{\pm}$.
Thus,~$Y$ is included in $\mu_Y(B^{\pm}_w(\lambda))$.
Let us def\/ine the Demazure polynomial $D^\pm_w[\lambda,Y;C]$ associated with a~monomial~$Y$, $w\in W$ and coef\/f\/icients
$C=(c(b))_{b\in B^\pm_w(\lambda)}$ ($c(b)\in \mathbb{Z}_{>0}$),
\begin{gather*}
D^{\pm}_w[\lambda,Y;C]:=\sum\limits_{b\in B^{\pm}_w(\lambda)}c(b)\mu_Y(b).
\end{gather*}
\end{defn}

\begin{rem}
In this paper, we only treat the case that the coef\/f\/icients $c(b)$ are equal to $1$ for all $b\in B^{-}_w(\lambda)$
(see Theorem~\ref{thm2}).
But when $G\neq {\rm SL}_{r+1}(\mathbb{C})$ (for example, $G={\rm Sp}_{2r}(\mathbb{C})$), we need to treat the case $c(b)$ is not
necessary equal to $1$ for some $b\in B^{-}_w(\lambda)$.
Therefore, we need non-trivial coef\/f\/icients $c(b)\in \mathbb{Z}_{>0}$ in Def\/inition~\ref{dempoly}.
\end{rem}

\section{Cluster algebras and generalized minors}

In this section, we shall review the notions of cluster algebras.
For all def\/initions in this section, see, e.g.,~\cite{A-F-Z, M-M-A}.

We set $[1,l]:=\{1,2,\dots,l\}$ and $[-1,-l]:=\{-1,-2,\dots,-l\}$ for $l\in \mathbb{Z}_{>0}$.
For $n,m\in \mathbb{Z}_{>0}$, let $x_1, \dots,x_n,x_{n+1}, \dots,x_{n+m}$ be variables and $\mathcal{P}$ be a~free
multiplicative abelian group generated by $x_{n+1},\dots,x_{n+m}$.
We set ${\mathbb Z}{\mathcal P}:={\mathbb Z}[x_{n+1}^{\pm1}, \dots,x_{n+m}^{\pm1}]$.
Let $K:=\{\frac{g}{h}  \,|\, g, h \in {\mathbb Z}{\mathcal P}$, $h\neq0 \}$ be the f\/ield of fractions of ${\mathbb Z}{\mathcal P}$,
and ${\mathcal F}:=K(x_{1}, \dots,x_{n})$ be the f\/ield of rational functions.

\subsection{Cluster algebras of geometric type}

\begin{defn}
We set~$n$-tuple of variables ${\bf x}=(x_1,\dots,x_n)$.
Let $\tilde{B}=(b_{ij})_{1\leq i\leq n,\; 1\leq j \leq n+m}$ be $n\times (n+m)$ integer matrix whose principal part
$B:=(b_{ij})_{1\leq i,j\leq n}$ is sign skew symmetric.
Then a~pair $\Sigma=({\bf x},\tilde{B})$ is called a~{\it seed}, {\bf x} a~cluster and $x_1, \dots, x_n$ cluster variables.
For a~seed $\Sigma=({\bf x},\tilde{B})$, principal part~$B$ of $\tilde{B}$ is called the {\it exchange matrix}.
\end{defn}

\begin{defn}
%\label{adc}
For a~seed $\Sigma=({\bf x}, \tilde{B}=(b_{ij}))$, an {\it adjacent cluster} in direction $k\in [1,n]$ is def\/ined~by
\begin{gather*}
{\bf x}_k = ({\bf x}\setminus \{x_k\})\cup \{x_k'\},
\end{gather*}
where $x_k'$ is the new cluster variable def\/ined by the {\it exchange relation}
\begin{gather*}
x_k x_k' = \prod\limits_{1\leq i \leq n+m,\; b_{ki}>0} x_i^{b_{ki}} +\prod\limits_{1\leq i \leq n+m,\; b_{ki}<0} x_i^{-b_{ki}}.
\end{gather*}
\end{defn}

\begin{defn}
Let $A=(a_{ij})$, $A'=(a_{ij}')$ be two matrices of the same size.
We say that $A'$ is obtained from~$A$ by the matrix mutation in direction~$k$, and denote $A'=\mu_k(A)$ if
\begin{gather*}
a_{ij}'=
\begin{cases}
-a_{ij} & \text{if} \quad i=k \quad \text{or} \quad j=k,
\\
a_{ij}+\dfrac{|a_{ik}|a_{kj}+a_{ik}|a_{kj}|}{2} & \text{otherwise}.
\end{cases}
\end{gather*}
For $A$, $A'$, if there exists a~f\/inite sequence $(k_1,\dots,k_s)$, $k_i\in[1,n]$, such that
$A'=\mu_{k_1}\cdots\mu_{k_s}(A)$, we say~$A$ is mutation equivalent to $A'$, and denote $A \cong A'$.
\end{defn}

Next proposition can be easily verif\/ied by the def\/inition of $\mu_k$:

\begin{prop}[\protect{\cite[Proposition~3.6]{M-M-A}}]\label{totally}
Let~$A$ be a~skew symmetrizable matrix.
Then any matrix that is mutation equivalent to~$A$ is sign skew symmetric.
\end{prop}

For a~seed $\Sigma=({\bf x},\tilde{B})$, we say that the seed $\Sigma'=({\bf x}',\tilde{B}')$ is adjacent
to~$\Sigma$ if ${\bf x}'$ is adjacent to ${\bf x}$ in direction~$k$ and $\tilde{B}'=\mu_k(\tilde{B})$.
Two seeds~$\Sigma$ and $\Sigma_0$ are mutation equivalent if one of them can be obtained from another seed by a~sequence
of pairwise adjacent seeds and we denote $\Sigma \sim\Sigma_0$.

Now let us def\/ine a~cluster algebra of geometric type.

\begin{defn}
Let $\tilde{B}$ be a~skew symmetrizable matrix, and $\Sigma=({\bf x},\tilde{B})$ a~seed.
We set ${\mathbb A}:={\mathbb Z}[x_{n+1},\dots,x_{n+m}]$.
The cluster algebra (of geometric type) ${\mathcal A}={\mathcal A}(\Sigma)$ over $\mathbb A$ associated with
seed~$\Sigma$ is def\/ined as the ${\mathbb A}$-subalgebra of ${\mathcal F}$ generated by all cluster variables in all
seeds which are mutation equivalent to~$\Sigma$.
\end{defn}

For a~seed~$\Sigma$, we def\/ine ${\mathbb Z}{\mathcal P}$-subalgebra ${\mathcal U}(\Sigma)$ of ${\mathcal F}$~by
\begin{gather*}
{\mathcal U}(\Sigma):={\mathbb Z}{\mathcal P}\big[{\bf x}^{\pm 1}\big] \cap {\mathbb Z}{\mathcal P}\big[{\bf x}_1^{\pm 1}\big]
\cap \dots \cap {\mathbb Z}{\mathcal P}\big[{\bf x}_n^{\pm 1}\big].
\end{gather*}
Here, ${\mathbb Z}{\mathcal P}[{\bf x}^{\pm 1}]$ is the Laurent polynomial ring in~${\bf x}$.

\begin{defn}
\label{upper}
Let $\Sigma_0=({\bf x},\tilde{B})$ be a~seed such that $\tilde{B}$ is skew symmetrizable.
We def\/ine an {\it upper cluster algebra} $\overline{{\mathcal A}}=\overline{{\mathcal A}}(\Sigma_0)$ as the intersection
of the subalgebras ${\mathcal U}(\Sigma)$ for all seeds $\Sigma \sim\Sigma_0$.
\end{defn}

Following the inclusion relation holds~\cite{A-F-Z}:
\begin{gather*}
{\mathcal A}(\Sigma) \subset \overline{{\mathcal A}}(\Sigma).
\end{gather*}

\subsection{Cluster algebras on double Bruhat cells of type A}

As in Section~\ref{DBCs}, let $G={\rm SL}_{r+1}(\mathbb{C})$ be the simple algebraic group of type $\text{A}_r$ and
$W=\mathfrak{S}_{r+1}$ be its Weyl group.
We set $u\in W$ and its reduced word ${\bf i}$ as in~\eqref{uset0} and~\eqref{iset0}:
\begin{gather}
\label{uset00}
u=\underbrace{s_1s_2\cdots s_r}_{\text{$1^{\rm st}$ cycle}} \underbrace{s_1\cdots s_{r-1}}_{\text{$2^{\rm nd}$ cycle}} \cdots
\underbrace{s_1\cdots s_{r-m+2}}_{\text{$(m-1)^{\rm th}$ cycle}} \underbrace{s_1\cdots s_{i_n}}_{\text{$m^{\rm th}$ cycle}},
\\
\label{iset00}
{\bf i}=(\underbrace{1,\dots,r}_{\text{$1^{\rm st}$ cycle}},\underbrace{1,\dots,(r-1)}_{\text{$2^{\rm nd}$ cycle}},\dots,\underbrace{1,\dots,(r-m+2)}_{\text{$(m-1)^{\rm th}$ cycle}}, \underbrace{1,\dots,i_n}_{\text{$m^{\rm th}$ cycle}}).
\end{gather}

We shall constitute the upper cluster algebra $\overline{{\mathcal A}}({\bf i})$ from ${\bf i}$.
Let $i_k$, $k\in[1,l(u)]$, be the~$k^{\rm th}$ index of ${\bf i}$ from the left.

At f\/irst, we def\/ine a~set $e({\bf i})$ as
\begin{gather*}
e({\bf i}):=[-1,-r]\cup \{k \,|\,  \text{there\ exist\ some}\ l>k\ \text{such\ that}\ i_k=i_l \}.
\end{gather*}

Next, let us def\/ine a~matrix $\tilde{B}=\tilde{B}({\bf i})$.

\begin{defn}
Let $\tilde{B}({\bf i})$ be an integer matrix with rows labeled by all the indices in $[-1,-r]\cup [1,l(u)]$ and
columns labeled by all the indices in $e({\bf i})$.
For $k\in[-1,-r]\cup [1,l(u)]$ and $l\in e({\bf i})$, an entry $b_{kl}$ of $\tilde{B}({\bf i})$ is determined as
follows
\begin{gather*}
b_{kl}=
\begin{cases}
-\sgn((k-l)\cdot i_p) & \text{if}\quad p=q,
\\
-\sgn((k-l)\cdot i_p\cdot a_{|i_k||i_l|}) & \text{if}\quad p<q \quad \text{and}\quad \sgn(i_p\cdot i_q)(k-l)(k^+-l^+)>0,
\\
0 & \text{otherwise}.
\end{cases}
\end{gather*}
\end{defn}

\begin{prop}[\protect{\cite[Proposition~2.6]{A-F-Z}}]\label{propss}
The matrix $\tilde{B}({\bf i})$ is skew symmetrizable.
\end{prop}

By Proposition~\ref{totally}, Def\/inition~\ref{upper} and Proposition~\ref{propss}, we can construct the upper
cluster algebra from $\tilde{B}({\bf i})$:

\begin{defn}
We denote this upper cluster algebra by $\overline{A}({\bf i})$.
\end{defn}

Now, we set $\bar{{\mathcal A}}({\bf i})_{\mathbb{C}}:=\bar{{\mathcal A}}({\bf i})\otimes \mathbb{C}$ and
${\mathcal F}_{\mathbb{C}}:={\mathcal F} \otimes \mathbb{C}$.
It is known that the coordinate ring $\mathbb{C}[G^{u,e}]$ of the double Bruhat cell is isomorphic to $\bar{{\mathcal
A}}({\bf i})_{\mathbb{C}}$ (Theorem~\ref{clmainthm}).
To describe this isomorphism explicitly, we need {\it generalized minors}.
For $k \in [1,l(u)]$, let $i_k$ be the~$k^{\rm th}$ index of ${\bf i}$ \eqref{iset00} from the left, and we suppose that
it belongs to the~$m'^{\rm th}$ cycle.
We set
\begin{gather}
\label{inc}
u_{\leq k}=u_{\leq k}({\bf i}):=\underbrace{s_1s_2\cdots s_r}_{\text{$1^{\rm st}$ cycle}} \underbrace{s_1\cdots s_{r-1}}_{\text{$2^{\rm nd}$ cycle}}
\cdots \underbrace{s_1\cdots s_{i_k}}_{\text{$m'^{\rm th}$ cycle}}.
\end{gather}
For $k \in [-1,-r]$, we set $u_{\leq k}:=e$ and $i_k:=k$.
In the case $G={\rm SL}_{r+1}(\mathbb{C})$, the generalized minors are nothing but the ordinary minors of a~matrix:
\begin{defn}[\cite{A-F-Z}]\label{Alem}
For $x \in G={\rm SL}_{r+1}(\mathbb{C})$ and $k \in [-1,-r]\cup[1,l(u)]$, we def\/ine the {\it generalized minor}
$\Delta(k;{\bf i})(x)$ as the minor of~$x$ whose rows (resp.\ columns) are labeled by the elements of the set
$u_{\leq k}([1,|i_k|])$ (resp.~$[1,|i_k|]$).
\end{defn}

Finally, we set
\begin{gather*}
F({\bf i}):=\{\Delta(k;{\bf i}) \,|\, k \in [-1,-r]\cup[1,l(u)] \}.
\end{gather*}
It is known that the set $F({\bf i})$ is an algebraically independent generating set for the f\/ield of rational
functions $\mathbb{C}(G^{u,e})$~\cite[Theorem 1.12]{F-Z}.
Then, we have the following theorem.

\begin{thm}[\protect{\cite[Theorem~2.10]{A-F-Z}}]\label{clmainthm}
The isomorphism of fields $\varphi:F_{\mathbb{C}} \rightarrow \mathbb{C}(G^{u,e})$
defined by $\varphi (x_k)=\Delta(k;{\bf i})$, $k \in [-1,-r]\cup [1,l(u)]$, restricts to an isomorphism of algebras
$\bar{{\mathcal A}}({\bf i})_{\mathbb{C}}\rightarrow \mathbb{C}[G^{u,e}]$.
\end{thm}

\section{Generalized minors and crystals}\label{gmc}

In the rest of the paper, we consider the case $G={\rm SL}_{r+1}(\mathbb{C})$, and let $u\in W$ and its reduced word
${\bf i}$ as in~\eqref{uset00} and~\eqref{iset00}:
\begin{gather}
\label{uset000}
u=s_1s_2\cdots s_rs_1\cdots s_{r-1}\cdots s_1\cdots s_{r-m+2} s_1\cdots s_{i_n},
\\
\label{iset000}
{\bf i}=(\underbrace{1,\dots,r}_{\text{$1^{\rm st}$~cycle}},\underbrace{1,\dots,(r-1)}_{\text{$2^{\rm nd}$
cycle}},\dots,\underbrace{1,\dots,(r-m+2)}_{\text{$(m-1)^{\rm th}$ cycle}}, \underbrace{1,\dots,i_n}_{\text{$m^{\rm th}$ cycle}}),
\end{gather}
that is, {\bf i} is the left factor of $(1,2,3,\dots,r,1,2,3,\dots,(r-1),\dots,1,2,1)$.
Let $i_k$ be the~$k^{\rm th}$ index of ${\bf i}$ from the left, and belong to $m'^{\rm th}$ cycle.
As we shall show in Lemma~\ref{gmlem}, we may assume $i_n=i_k$.

By Theorem~\ref{clmainthm}, we can regard $\mathbb{C}[G^{u,e}]$ as an upper cluster algebra and
$\{\Delta(k;{\bf i})\}$ as its cluster variables belonging to a~given initial seed.
Each $\Delta(k;{\bf i})$ is a~regular function on $G^{u,e}$.
On the other hand, by Proposition~\ref{gprime} (resp.\ Theorem~\ref{fp2}), we can consider $\Delta(k;{\bf i})$
as a~function on $H\times (\mathbb{C}^{\times})^{l(u)}$
(resp.\ $(\mathbb{C}^{\times})^{l(u)}$).
Then we change the variables of $\{\Delta(k;{\bf i})\}$ as follows:
\begin{defn}
%\label{gendef}
For $a\in H$ and ${\bf t}$, $\tau\in (\mathbb{C}^{\times})^{l(u)}$ we set
\begin{gather*}
\Delta^G(k;{\bf i})(a,{\bf t}):=\big(\Delta(k;{\bf i})\circ \bar{x}^G_{{\bf i}}\big)(a,{\bf t}),
\\
\Delta^L(k;{\bf i})(\tau):=\big(\Delta(k;{\bf i})\circ x^L_{{\bf i}}\big)(\tau),
\end{gather*}
where ${\bf t}=(t_1,\dots,t_{l(u)})$, $\tau=(\tau_1,\dots,\tau_{l(u)})$.
\end{defn}
We will describe the function $\Delta^L(k;{\bf i})(\tau)$ by using monomial realizations of Demazure crystals.

\subsection[Generalized minor $\Delta^G(k;{\bf i})(a,{\bf t})$]{Generalized minor $\boldsymbol{\Delta^G(k;{\bf i})(a,{\bf t})}$}

In this subsection, we shall prove that $\Delta^G(k;{\bf i})(a,{\bf t})$ is immediately obtained from
$\Delta^L(k;{\bf i})$:
\begin{prop}
\label{gprop}
We set $d:=i_k$.
For $a=\diag(a_1,\dots,a_{r+1})\in H$,
\begin{gather*}
\Delta^G(k;{\bf i})(a,{\bf t})=a_{m'+1}\cdots a_{m'+d}\Delta^L(k;{\bf i})({\bf t}).
\end{gather*}
\end{prop}
This proposition follows from the following lemma:
\begin{lem}
\label{gllem}
In the above setting, $\Delta^G(k;{\bf i})(a,{\bf t})$ $($resp.~$\Delta^L(k;{\bf i})(\tau))$ is given as a~minor whose row are labeled
by the set $\{m'+1,\dots,m'+d\}$ and column are labeled
by the set $\{1,\dots,d\}$ of the matrix
\begin{gather*}
a\underbrace{x_{-1}(t_1)x_{-2}(t_2)\cdots x_{-r}(t_r)}_{\text{\rm $1^{\rm st}$ cycle}}\underbrace{x_{-1}(t_{l_1+1})\cdots
x_{-(r-1)}(t_{l_1+r-1})}_{\text{\rm $2^{\rm nd}$ cycle}}\cdots
\\
\qquad{} \times\underbrace{x_{-1}(t_{l_{m-1}+1}) \cdots x_{-i_n}(t_{l_{m-1}+i_n})}_{\text{\rm $m^{\rm th}$ cycle}},
\end{gather*}
resp.\
\begin{gather}
\underbrace{x_{-1}(\tau_1)x_{-2}(\tau_2)\cdots x_{-r}(\tau_r)}_{\text{\rm $1^{\rm st}$ cycle}}\underbrace{x_{-1}(\tau_{l_1+1})\cdots
x_{-(r-1)}(\tau_{l_1+r-1})}_{\text{\rm $2^{\rm nd}$ cycle}}\cdots
\nonumber
\\
\qquad {}\times \underbrace{x_{-1}(\tau_{l_{m-1}+1})\cdots x_{-i_n}(\tau_{l_{m-1}+i_n})}_{\text{\rm $m^{\rm th}$ cycle}}.
\label{tmatrix}
\end{gather}
\end{lem}

\begin{proof} Let us prove this lemma for $\Delta^L(k;{\bf i})(\tau)$ since the case for
$\Delta^G(k;{\bf i})(a,{\bf t})$ is proven similarly.
By the def\/inition \eqref{inc} of $u_{\leq k}$ and $i_k=d$, we have
\begin{gather*}
u_{\leq k}[1,|i_k|]=u_{\leq k}\{1,\dots,d\} = \underbrace{s_1\cdots s_r}_{\text{$1^{\rm st}$ cycle}}\cdots \underbrace{s_1\cdots
s_{r-m'+2}}_{\text{$(m'-1)^{\rm th}$ cycle}}\underbrace{s_1\cdots s_d}_{\text{$m'^{\rm th}$ cycle}}\{1,\dots,d\}
\\
\phantom{u_{\leq k}[1,|i_k|]}
 = \{m'+1,\dots,m'+d\}.
\end{gather*}
Hence, by Theorem~\ref{fp2} and Def\/inition~\ref{Alem},
$\Delta^L(k;{\bf i})(\tau)$ is given as a~minor whose row (resp.\ column) are labeled
by the set $\{m'+1,\dots,m'+d\}$ (resp.\ $\{1,\dots,d\}$) of the matrix
\begin{gather*}
\underbrace{x_{-1}(\tau_1)x_{-2}(\tau_2)\cdots x_{-r}(\tau_r)}_{\text{$1^{\rm st}$ cycle}}\underbrace{x_{-1}(\tau_{l_1+1})\cdots
x_{-(r-1)}(\tau_{l_1+r-1})}_{\text{$2^{\rm nd}$ cycle}}\cdots
\\
\qquad{}\times \underbrace{x_{-1}(\tau_{l_{m-1}+1})\cdots x_{-i_n}(\tau_{l_{m-1}+i_n})}_{\text{$m^{\rm th}$ cycle}},
\end{gather*}
which implies the desired result.
\end{proof}

In the rest of the paper, we will treat $\Delta^L(k;{\bf i})(\tau)$ only by Proposition~\ref{gprop}.

\subsection[Generalized minor $\Delta^L(k;{\bf i})(\tau)$]{Generalized minor $\boldsymbol{\Delta^L(k;{\bf i})(\tau)}$}\label{Mainsec}

\begin{lem}
\label{gmlem}
Let~$u$ and ${\bf i}=(1,\dots, i_n)$ be as in the form \eqref{uset000} and \eqref{iset000} respectively.
For $i_{n+1}\in [1,r]$ and $u':=us_{i_{n+1}}\in W$ ($l(u')>l(u)$)
we set the reduced word for $u'$ as
${\bf i}':=(1,\dots, i_n,i_{n+1})$.
We also set $\tau=(\tau_1,\dots,\tau_n)$ and $\tau'=(\tau_1,\dots,\tau_n,\tau_{n+1})$.
For an integer~$k$, $1\leq k\leq n$, if $d:=i_k\neq i_{n+1}$, then $\Delta^L(k;{\bf i}')(\tau')$ does not depend on
$\tau_{n+1}$, so we can regard it as a~function on $H\times (\mathbb{C}^{\times})^{n}$.
Furthermore, we have
\begin{gather*}
%\label{Lomit}
\Delta^L(k;{\bf i})(\tau)=\Delta^L(k;{\bf i}')(\tau').
\end{gather*}
\end{lem}

\begin{proof} We denote the matrix~\eqref{tmatrix} by $T=(T_{i,j})_{i,j\in [1,r+1]}$:
\begin{gather*}
T=(T_{i,j})_{i,j\in [1,r+1]}:=x_{-1}(\tau_1)x_{-2}(\tau_2)\cdots x_{-i_n}(\tau_{n}).
\end{gather*}
We also def\/ine the submatrix $T_u$ of the matrix~$T$ whose rows (resp.\ columns) are labeled
by the set $\{m'+1,\dots,m'+d\}$ (resp.\ $\{1,\dots,d\}$), that is
\begin{gather*}
T_u:=
\begin{pmatrix}
T_{m'+1, 1} & \dots & T_{m'+1, d}
\\
\vdots & \vdots & \vdots
\\
T_{m'+d, 1} & \dots & T_{m'+d, d}
\end{pmatrix}
.
\end{gather*}
By Lemma~\ref{gllem}, $\Delta^L(k;{\bf i})(\tau)$ is the determinant of $T_u$.

Similarly, we def\/ine the submatrix $T'_{u'}$ of the matrix $T':=T\cdot x_{-i_{n+1}}(\tau_{n+1})$ whose rows (resp.\ columns) are labeled
by the set $\{m'+1,\dots,m'+d\}$ (resp.\ $\{1,\dots,d\}$).
By Lemma~\ref{gllem}, $\Delta^L(k;{\bf i}')(\tau')$ is the determinant of $T'_{u'}$.

Using the explicit form of $x_{-i}(t)$ in~\eqref{transmat}, we have the following relation: If $i_{n+1}<i_k=d$, then
\begin{gather*}
T'_{u'} =
\begin{pmatrix}
T_{m'+1,1} & \dots & \tau_{n+1}^{-1} T_{m'+1, i_{n+1}}+T_{m'+1, i_{n+1}+1} & \tau_{n+1} T_{m'+1,i_{n+1}+1} & \dots &
T_{m'+1,d}
\\
\vdots & \vdots & \vdots & \vdots & \vdots & \vdots
\\
T_{m'+d, 1} & \dots & \tau_{n+1}^{-1} T_{m'+1, i_{n+1}}+T_{m'+1, i_{n+1}+1} & \tau_{n+1}T_{m'+1, i_{n+1}+1} & \dots &
T_{m'+d,d}
\end{pmatrix}
.
\end{gather*}
If $i_{n+1}>i_k=d$, then $T'_{u'}=T_u$.
Therefore, it is clear that for both cases, $\det T'_{u'}=\det T_u$, which means
$\Delta^L(k;{\bf i}')(\tau')=\Delta^L(k;{\bf i})(\tau)$.
\end{proof}

By this lemma, when we calculate $\Delta^L(k;{\bf i})(\tau)$, we may assume that $i_n=i_k$ without loss of
generality.

We set $l_i=\sum\limits^{i}_{k=1}(r-k+1) (1\leq i\leq m-1)$, and change the variables $\{Y_{m,j}\}$ to
$\{\tau_{l_{m}+j}\}$ as in Section~\ref{monoreal}.

\begin{rem}
For
\begin{gather*}
{\bf i}=(\underbrace{1,\dots,r}_{\text{$1^{\rm st}$ cycle}}, \underbrace{1,\dots,(r-1)}_{\text{$2^{\rm nd}$  cycle}},\dots,\underbrace{1,\dots,(r-m+2)}_{\text{$(m-1)^{\rm th}$ cycle}},\underbrace{1,\dots,i_n}_{\text{$m^{\rm th}$ cycle}}),
\end{gather*}
$(l_i+j)^{\rm th}$ index of {\bf i} from the left is~$j$ which belongs to $(i+1)^{\rm th}$ cycle $(1\leq j\leq r-i-1)$.
\end{rem}
The following theorem is our main result.
We describe $\Delta^L(k;{\bf i})(\tau)$ as a~Demazure polyno\-mial~$D^-_w$ in Def\/inition~\ref{dempoly}.

\begin{thm}
\label{thm2}
In the above setting, we set $d:=i_k=i_n$ and suppose that $i_k$ belongs to $m'^{\rm th}$ cycle in~${\bf i}$
\eqref{iset000}.
We also set
\begin{gather*}
%\label{yy}
Y:=\frac{1}{\tau_{l_{m-1}+d}\tau_{l_{m-2}+d}\cdots\tau_{l_{m'}+d}}.
\end{gather*}
Then we have
\begin{gather*}
\Delta^L(k;{\bf i})(\tau)  = \sum\limits_{x\in B^-_{u_{\leq k}}((m'-m)\Lambda_d)}\mu_Y(x)
 = D^-_{u_{\leq k}}[(m'-m)\Lambda_d,Y;{\bf 1}],
\end{gather*}
where $\mu_Y(x)$ is an embedding of $x\in B^-_{u_{\leq k}}((m'-m)\Lambda_d)$ as in Example~{\rm \ref{emb}}, and ${\bf 1}$
means all $c(b)\equiv1$ $($Definition~{\rm \ref{dempoly})}.
\end{thm}

We shall prove this theorem in Section~\ref{prsec}.
In particular, we can explicitly write down $\Delta^L(k;{\bf i})(\tau)$ in the case $i_k=1$.

\begin{thm}
\label{thm3}
If $i_k=1$ and $i_k$ is in the $m'^{\rm th}$ cycle in ${\bf i}$ \eqref{iset000}, then we have
\begin{gather*}
\Delta^L(k;{\bf i})(\tau)
= \sum\limits_{0\leq j_1<\dots<j_{m'}\leq m-1} \prod\limits^{j_1-1}_{i=0} \frac{1}{\tau_{l_{m-1-i}+1}}
\prod\limits^{j_2-1}_{i=j_1+1}\frac{\tau_{l_{m-1-i}+1}}{\tau_{l_{m-1-i}+2}}\cdots
\prod\limits^{m-1}_{i=j_{m'}+1}\frac{\tau_{l_{m-1-i}+m'}}{\tau_{l_{m-1-i}+m'+1}}
\\
\phantom{\Delta^L(k;{\bf i})(\tau)}
=D^-_{u_{\leq k}}[(m'-m)\Lambda_1, \frac{1}{\tau_{l_{m-1}+1}\tau_{l_{m-2}+1}\cdots \tau_{l_{m'+1}+1}
\tau_{l_{m'}+1}};{\bf 1}].
\end{gather*}
\end{thm}

We will prove this theorem in Section~\ref{prthm3}.

\section{The proof of Theorems~\ref{thm2} and~\ref{thm3}}\label{prsec}

In this section, we shall give the proof of Theorems~\ref{thm2} and~\ref{thm3}.
In Sections~\ref{prthm2start}--\ref{prthm2fin}, we will prove Theorem~\ref{thm2}.
In Section~\ref{prthm3}, we prove Theorem~\ref{thm3}.
We use the same notation as in Section~\ref{Mainsec}: $l_i=\sum\limits^{i}_{k=1}(r-k+1)$, $1\leq i\leq m-1$.

\subsection[The set $X_d(m,m')$ of paths]{The set $\boldsymbol{X_d(m,m')}$ of paths}\label{prthm2start}

In this subsection, we shall introduce a~set $X_d(m,m')$ of ``paths'' which correspond to the terms of
$\Delta^L(k;{\bf i})(\tau)$.
Let~$m$, $m'$ and~$d$ be the positive integers as in Section~\ref{Mainsec}.
\begin{defn}
Let us def\/ine the directed graph $(V_d,E_d)$ as follows. We def\/ine the set $V_d=V_d(m,m')$ of vertices as
\begin{gather*}
V_d(m,m'):=\big\{\big(m-s;a^{(s)}_1,a^{(s)}_2,\dots,a^{(s)}_d\big) \,\big|\,
\\
\hphantom{V_d(m,m'):=\big\{}
0\leq s\leq m,\ a^{(s)}_i\in \mathbb{Z},\ 1\leq
a^{(s)}_1<a^{(s)}_2<\dots<a^{(s)}_d\leq d+s \big\}.
\end{gather*}
Note that $a^{(0)}_1=1$, $a^{(0)}_2=2,\dots,a^{(0)}_d=d$ by $1\leq a^{(0)}_1<a^{(0)}_2<\dots<a^{(0)}_d\leq d$.
And we def\/ine the set $E_d=E_d(m,m')$ of directed edges as
\begin{gather*}
E_d(m,m')
:=\big\{\big(m-s;a^{(s)}_1,\dots,a^{(s)}_d\big)
\rightarrow \big(m-s-1;a^{(s+1)}_1,\dots,a^{(s+1)}_d\big) \;\big|\; \\
\hphantom{E_d(m,m') :=\big\{}{}
1\leq s\leq m,\ 1\leq i\leq d,\ a^{(s+1)}_i
=a^{(s)}_i\ \text{or}\ a^{(s)}_i+1\big\}.
\end{gather*}
\end{defn}

\begin{defn}
\label{pathdef}
Let $X_d(m,m')$ be the set of directed paths from $(m;1,\dots,d)$ to $(0;m'+1$, $m'+2,\dots,m'+d)$ in $(V_d,E_d)$.
In other word, any path $p\in X_d(m,m')$
\begin{gather*}
p=\big(m;a^{(0)}_1,\dots,a^{(0)}_d\big) \rightarrow\big(m-1;a^{(1)}_1,\dots,a^{(1)}_d\big)\rightarrow\big(m-2;a^{(2)}_1,
\dots,a^{(2)}_d\big)\rightarrow\cdots
\\
\hphantom{p=}{} \rightarrow\big(1;a^{(m-1)}_1,\dots,a^{(m-1)}_d\big) \rightarrow\big(0;a^{(m)}_1,\dots,a^{(m)}_d\big)
\end{gather*}
are characterized by the following conditions
\begin{enumerate}\itemsep=0pt
\item[(i)] $a^{(s)}_{i}\in \mathbb{Z}_{\geq 1}$,
\item[(ii)] $a^{(s)}_{1}<a^{(s)}_{2}<\dots<a^{(s)}_{d}$,
\item[(iii)] $a^{(s+1)}_{i}=a^{(s)}_{i}$ or $a^{(s)}_{i}+1$,
\item[(iv)] $a^{(0)}_i=i$, $a^{(m)}_i=m'+i$.
\end{enumerate}
\end{defn}

\begin{rem}
\label{pathrem}
By Def\/inition~\ref{pathdef}(ii), (iii) and (iv), we have
\begin{gather*}
a^{(s)}_i\leq a^{(s)}_d \leq a^{(m)}_d=m'+d
\end{gather*}
for any $1\leq i\leq d$ and $0\leq s\leq m$.
\end{rem}
Let us def\/ine a~Laurent monomial associated with each path in $X_d(m,m')$.
\begin{defn}
\label{labeldef}
Let $p\in X_d(m,m')$ be a~path
\begin{gather*}
p=\big(m;a^{(0)}_1,\dots,a^{(0)}_d\big)\rightarrow\big(m-1;a^{(1)}_1,
\dots,a^{(1)}_d\big)\rightarrow\big(m-2;a^{(2)}_1,\dots,a^{(2)}_d\big)\rightarrow\cdots
\\
\hphantom{p=}{}
\rightarrow\big(1;a^{(m-1)}_1,\dots, a^{(m-1)}_d\big)\rightarrow\big(0;a^{(m)}_1,\dots,a^{(m)}_d\big).
\end{gather*}
\begin{enumerate}\itemsep=0pt
\item[(i)] For each $0\leq s\leq m$, we def\/ine the {\it label of the edge}
$\big(m-s;a^{(s)}_1,a^{(s)}_2,\dots,a^{(s)}_d\big)\rightarrow\big(m-s-1;$ $a^{(s+1)}_1,a^{(s+1)}_2,\dots,a^{(s+1)}_d\big)$ as the Laurent
monomial
\begin{gather}
\label{label}
\prod\limits^d_{i=1}\frac{\tau_{l_{m-s-1}+a^{(s+1)}_i-1}}{\tau_{l_{m-s-1}+a^{(s)}_i}}.
\end{gather}
\item[(ii)] And we def\/ine the {\it label} $Q(p)$ {\it of the path}~$p$ as the product of them
\begin{gather*}
%\label{alllabels}
Q(p):=\prod\limits_{s=0}^{m-1}\left(\prod\limits^d_{i=1}\frac{\tau_{l_{m-s-1}+a^{(s+1)}_i-1}}{\tau_{l_{m-s-1}+a^{(s)}_i}}
\right).
\end{gather*}
\end{enumerate}
\end{defn}

\begin{ex}
\label{pathex}
Let $m=3$, $m'=2$, $d=2$.
We can describe $X_2(3,2)$ and its labels as follows
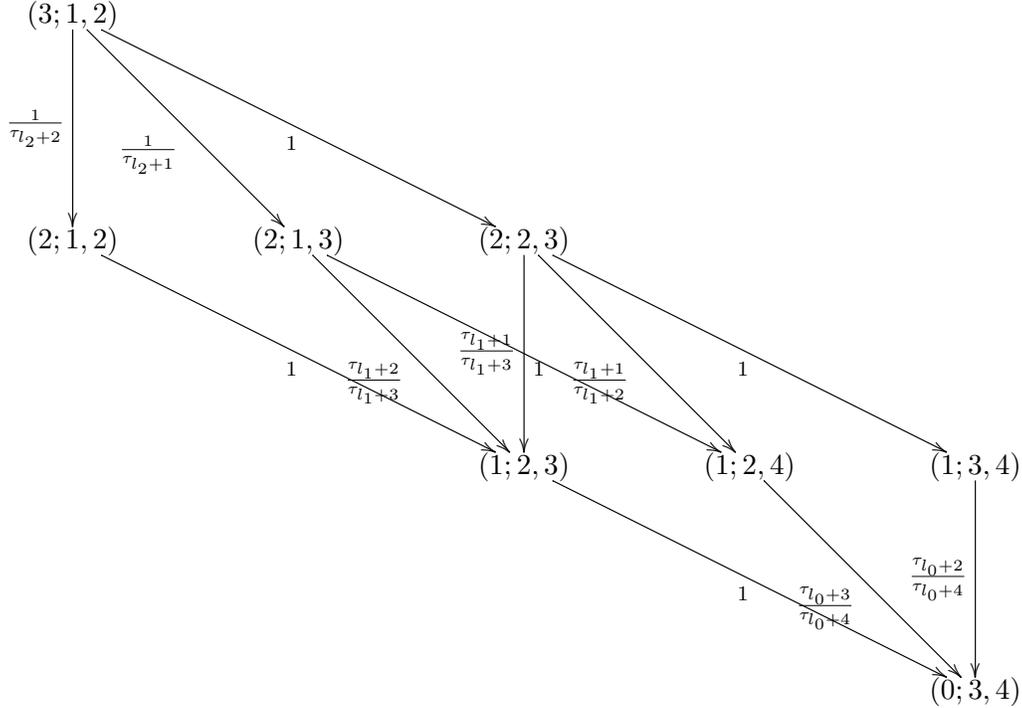
\begin{gather*}
\begin{xy}
(0,90) *{(3;1,2)}="3;1,2",
(0,60)*{(2;1,2)}="2;1,2",
(30,60)*{(2;1,3)}="2;1,3",
(60,60)*{(2;2,3)}="2;2,3",
(60,30)*{(1;2,3)}="1;2,3",
(90,30)*{(1;2,4)}="1;2,4",
(120,30)*{(1;3,4)}="1;3,4",
(120,0)*{(0;3,4)}="0;3,4",
\ar@{->} "3;1,2";"2;1,2"_{\frac{1}{\tau_{l_{2}+2}}}
\ar@{->} "3;1,2";"2;1,3"_{\frac{1}{\tau_{l_2+1}}}
\ar@{->} "3;1,2";"2;2,3"_1
\ar@{->} "2;1,2";"1;2,3"_1
\ar@{->} "2;1,3";"1;2,3"_{\frac{\tau_{l_1+2}}{\tau_{l_1+3}}}
\ar@{->} "2;1,3";"1;2,4"_{\qquad \quad 1}                      %\qquad \quad
\ar@{->} "2;2,3";"1;2,3"_{\frac{\tau_{l_1+1}}{\tau_{l_1+3}}}
\ar@{->} "2;2,3";"1;2,4"_{\frac{\tau_{l_1+1}}{\tau_{l_1+2}}}
\ar@{->} "2;2,3";"1;3,4"_1
\ar@{->} "1;2,3";"0;3,4"_1
\ar@{->} "1;2,4";"0;3,4"_{\frac{\tau_{l_0+3}}{\tau_{l_0+4}}}
\ar@{->} "1;3,4";"0;3,4"_{\frac{\tau_{l_0+2}}{\tau_{l_0+4}}}
\end{xy}
\end{gather*}
Each edge has the label written on the left side of it.
The paths in $X_2(3,2)$ are as follows
\begin{gather*}
p_1=(3;1,2)\rightarrow(2;1,2)\rightarrow(1;2,3)\rightarrow(0;3,4),\\
p_2=(3;1,2)\rightarrow(2;1,3)\rightarrow(1;2,3)\rightarrow(0;3,4),\\
p_3=(3;1,2)\rightarrow(2;1,3)\rightarrow(1;2,4)\rightarrow(0;3,4),\\
p_4=(3;1,2)\rightarrow(2;2,3)\rightarrow(1;2,3)\rightarrow(0;3,4),\\
p_5=(3;1,2)\rightarrow(2;2,3)\rightarrow(1;2,4)\rightarrow(0;3,4),\\
p_6=(3;1,2)\rightarrow(2;2,3)\rightarrow(1;3,4)\rightarrow(0;3,4).
\end{gather*}

We have
\begin{gather*}
Q(p_1)=\frac{1}{\tau_{l_{2}+2}},
\qquad
Q(p_2)=\frac{\tau_{l_1+2}}{\tau_{l_2+1}\tau_{l_1+3}},
\qquad
Q(p_3)=\frac{\tau_{l_0+3}}{\tau_{l_2+1}\tau_{l_0+4}},
\\
Q(p_4)=\frac{\tau_{l_1+1}}{\tau_{l_1+3}},
\qquad
Q(p_5)=\frac{\tau_{l_1+1}\tau_{l_0+3}}{\tau_{l_1+2}\tau_{l_0+4}},
\qquad
Q(p_6)=\frac{\tau_{l_0+2}}{\tau_{l_0+4}}.
\end{gather*}
\end{ex}

\begin{defn}
\label{iseq}
For each path $p\in X_d(m,m')$
\begin{gather*}
p=\big(m;a^{(0)}_1,\dots,a^{(0)}_d\big)\rightarrow\big(m-1;a^{(1)}_1,\dots,a^{(1)}_d\big)\rightarrow\big(m-2;a^{(2)}_1,\dots,a^{(2)}_d\big)\rightarrow\cdots
\\
\hphantom{p=}{}
\rightarrow\big(1;a^{(m-1)}_1,\dots,a^{(m-1)}_d\big)\rightarrow\big(0;a^{(m)}_1,\dots,a^{(m)}_d\big)
\end{gather*}
and $i\in \{1,\dots,d\}$, we call the following sequence
\begin{gather*}
a^{(0)}_i\rightarrow a^{(1)}_i\rightarrow a^{(2)}_i\rightarrow\cdots\rightarrow a^{(m)}_i
\end{gather*}
an {\it i-sequence} of~$p$.
\end{defn}

\begin{ex}
In the setting of Example~\ref{pathex}, $p_1:=(3;1,2)\rightarrow(2;1,2)\rightarrow(1;2,3)\rightarrow(0;3,4)$.
Then, $1$-sequence of $p_1$ is
\begin{gather*}
1\rightarrow1\rightarrow2\rightarrow3
\end{gather*}
and $2$-sequence of $p_1$ is
\begin{gather*}
2\rightarrow2\rightarrow3\rightarrow4.
\end{gather*}
\end{ex}

\subsection[One-to-one correspondence between paths in $X_d(m,m')$ and terms of $\Delta^L(k;{\bf i})(\tau)$]{One-to-one correspondence
between paths in $\boldsymbol{X_d(m,m')}$\\ and terms of $\boldsymbol{\Delta^L(k;{\bf i})(\tau)}$}

\begin{prop}
\label{pathlem}
We use the setting \eqref{uset000}, \eqref{iset000} and the notations in Section~{\rm \ref{gmc}}.
Then, we have the following
\begin{gather*}
\Delta^L(k;{\bf i})(\tau)=\sum\limits_{p\in X_d(m,m')} Q(p).
\end{gather*}
\end{prop}

To prove this proposition, we need the following preparations.
For~$s$ ($0\leq s \leq m-1$), we def\/ine a~matrix $x^{(s)}(\tau)$ as
\begin{gather*}
x^{(s)}(\tau):=\underbrace{x_{-1}(\tau_1)x_{-2}(\tau_2)\cdots x_{-r}(\tau_r)}_{\text{$1^{\rm st}$ cycle}}\underbrace{x_{-1}(\tau_{l_1+1})\cdots x_{-(r-1)}(\tau_{l_1+r-1})}_{\text{$2^{\rm nd}$ cycle}}\cdots
\\
\hphantom{x^{(s)}(\tau):=}{}\times\underbrace{x_{-1}(\tau_{l_{s-1}+1})\cdots x_{-(r-s+1)}(\tau_{l_{s-1}+r-s+1})}_{\text{$s^{\rm th}$ cycle}},
\end{gather*}
where we understand $x^{(0)}(\tau)$ means identity matrix.
We also def\/ine
\begin{gather*}
x^{(m)}(\tau):=\underbrace{x_{-1}(\tau_1)x_{-2}(\tau_2)\cdots x_{-r}(\tau_r)}_{\text{$1^{\rm st}$ cycle}}\underbrace{x_{-1}(\tau_{l_1+1})\cdots x_{-(r-1)}(\tau_{l_1+r-1})}_{\text{$2^{\rm nd}$ cycle}}\cdots
\\
 \hphantom{x^{(m)}(\tau):=}{} \times \underbrace{x_{-1}(\tau_{l_{m-1}+1})\cdots x_{-i_n}(\tau_{l_{m-1}+i_n})}_{\text{$m^{\rm th}$ cycle}},
\end{gather*}
which is equal to the matrix in~\eqref{tmatrix}.

For $x^{(s)}(\tau)=(x^{(s)}_{i,j})_{i,j\in [1,r+1]}$, we def\/ine the~$d$ dimensional column vector $D(s;p)$
\begin{gather*}
D(s;p):=
\begin{pmatrix}
x^{(s)}_{m'+1,p}
\\
\vdots
\\
x^{(s)}_{m'+d,p}
\end{pmatrix}
\in \mathbb{C}^d.
\end{gather*}
Note that, by the explicit form of $x_{-i}(\tau)$ in~\eqref{transmat}, multiplying $x_{-i}(\tau)$ from the right gives
an elementary transformation of a~matrix.
Therefore, we get
\begin{gather}
\label{xsp1}
D(s;p)=
\begin{cases}
\dfrac{1}{\tau_{l_{s-1}+p}}D(s-1;p)+D(s-1;p+1) & \text{if} \quad  p=1, \vspace{1mm}\\
\dfrac{\tau_{l_{s-1}+p-1}}{\tau_{l_{s-1}+p}}D(s-1;p)+D(s-1;p+1) & \text{if} \quad  p>1, \end{cases}
\end{gather}
for $1\leq s\leq m$.
For $0\leq s\leq m$ and $1\leq i_1<\dots<i_d \leq r$, we set
\begin{gather}
\label{sid}
(s;i_1,i_2,\dots,i_d):=\det(D(s;i_1),D(s;i_2),\dots,D(s;i_d)),
\end{gather}
which coincides with the notation for a~vertex in $X_d(m,m')$ since later we identify each vertex with the minor above.

\begin{proof}[Proof of Proposition~\ref{pathlem}]
We shall prove the proposition in the following three steps.

Step~1.
$\Delta^L(k;{\bf i})(\tau)=\det(D(m;1),D(m;2),\dots,D(m;d))$.

It is followed from Lemma~\ref{gllem} that $\Delta^L(k;{\bf i})(\tau)$ is given as a~minor whose row (resp.\ column) are labeled
by the set $\{m'+1,\dots,m'+d\}$ (resp.\ $\{1,\dots,d\}$) of the matrix $x^{(m)}(\tau)$.
Thus, by using above notation, we have
$\Delta^L(k;{\bf i})(\tau)=\det(^t(D(m;1),D(m;2),\dots,D(m;d)))=\det(D(m;1),D(m;2),\dots,D(m;d))$.

Step~2.
Calculation of $\det(D(m;1),D(m;2),\dots,D(m;d))$ and labeled
graph.

By using~\eqref{xsp1}, let us calculate $\Delta^L(k;{\bf i})(\tau)=\det(D(m;1),D(m;2),\dots,D(m;d))$
explicitly.
We have
\begin{gather}
  (m;1,2,\dots,d)=\det(D(m;1),D(m;2),\dots,D(m;d))
\nonumber
\\
\hphantom{(m;1,2,\dots,d)=}{}
 =\det\left(\frac{1}{\tau_{l_{m-1}+1}}D(m-1;1)+D(m-1;2), \dots,\right.
\nonumber
\\
\left.\hphantom{(m;1,2,\dots,d)=\det \ \ \ \ }{}
  \frac{\tau_{l_{m-1}+d-1}}{\tau_{l_{m-1}+d}}D(m-1;d)+D(m-1;d+1)\right),
\label{ceq}
\end{gather}
which implies that in the notation~\eqref{sid}, $(m;1,2,\dots,d)$ is a~linear combination of
$\{(m-1$; $a^{(1)}_1,a^{(1)}_2,\dots,a^{(1)}_d) \,|\, 1\leq a^{(1)}_1<a^{(1)}_2<\dots<a^{(1)}_d \leq d+1$, $a^{(1)}_i=i$ or $i+1\}$.
By~\eqref{ceq}, the coef\/f\/icient of each $(m-1;a^{(1)}_1,a^{(1)}_2,\dots,a^{(1)}_d)$ is
\begin{gather*}
\prod\limits^d_{i=1}\frac{\tau_{l_{m-1}+a^{(1)}_i-1}}{\tau_{l_{m-1}+i}},
\end{gather*}
which coincides with the label of the edge connecting $(m;1,\dots,d)$ and $(m-1;a^{(1)}_1,a^{(1)}_2,\dots,a^{(1)}_d)$
def\/ined in Def\/inition~\ref{labeldef}, see~\eqref{label}.

Using the formula~\eqref{xsp1} again, we see that each $(m-1;a^{(1)}_1,a^{(1)}_2,\dots,a^{(1)}_d)$ is a~linear
combination of $\{(m-2;a^{(2)}_1,a^{(2)}_2,\dots,a^{(2)}_d) \,|\, 1\leq a^{(2)}_1<a^{(2)}_2<\dots<a^{(2)}_d \leq d+2$,
$a^{(2)}_i=a^{(1)}_i$ or $a^{(1)}_i+1\}$ in the same way as~\eqref{ceq}.
Then, the coef\/f\/icient of $(m-2;a^{(2)}_1,a^{(2)}_2,\dots,a^{(2)}_d)$ is
\begin{gather*}
\prod\limits^d_{i=1}\frac{\tau_{l_{m-2}+a^{(2)}_i-1}}{\tau_{l_{m-2}+a^{(1)}_i}},
\end{gather*}
which coincides with the label of the edge connecting $(m-1;a^{(1)}_1,a^{(1)}_2,\dots,a^{(1)}_d)$ and
$(m-2$; $a^{(2)}_1,a^{(2)}_2,\dots,a^{(2)}_d)$.
Thus $(m;1,2,\dots,d)$ is a~linear combination of $(m-2;a^{(2)}_1,a^{(2)}_2,\dots,a^{(2)}_d)$,
$1\leq a^{(2)}_1<a^{(2)}_2<\dots<a^{(2)}_d \leq d+2$, whose coef\/f\/icient is
\begin{gather*}
\sum\limits_{a^{(1)}_1,\dots,a^{(1)}_d} \left(\prod\limits^d_{i=1}\frac{\tau_{l_{m-1}+a^{(1)}_i-1}}{\tau_{l_{m-1}+i}}
\times\prod\limits^d_{i=1}\frac{\tau_{l_{m-2} +a^{(2)}_i-1}}{\tau_{l_{m-2}+a^{(1)}_i}}\right),
\end{gather*}
where the index $a^{(1)}_i$ ($1\leq i\leq d$) runs over $\{(a^{(1)}_i:i\in\{1,\dots,d\}) \,|\, a^{(2)}_i=a^{(1)}_i$ or $a^{(1)}_i+1$, $a^{(1)}_i=i$
or $i+1$, $a^{(1)}_i< a^{(1)}_{i+1}\}$.

Repeating this argument, we see that $(m;1,2,\dots,d)$ is a~linear combination of
$(0;a^{(m)}_1,a^{(m)}_2$, $\dots,a^{(m)}_d)$, $1\leq a^{(m)}_1<a^{(m)}_2<\dots<a^{(m)}_d \leq m+d$, whose coef\/f\/icient is
\begin{gather}
\label{lincom1}
\sum\limits_{a^{(j)}_i,\ i\in \{1,\dots,d\},\ 1\leq j\leq m-1}
\left(\prod\limits^d_{i=1}\frac{\tau_{l_{m-1}+a^{(1)}_i-1}}{\tau_{l_{m-1}+i}}\times\prod\limits^d_{i=1}
\frac{\tau_{l_{m-2}+a^{(2)}_i-1}}{\tau_{l_{m-2}+a^{(1)}_i}} \times\dots\times
\prod\limits^d_{i=1}\frac{\tau_{l_{0}+a^{(m)}_i-1}}{\tau_{l_{0}+a^{(m-1)}_i}} \right),
\end{gather}
where the index $a^{(j)}_i$ runs over $\{(a^{(j)}_i:i\in \{1,\dots,d\}$, $1\leq j\leq m-1) \,|\, a^{(j+1)}_i=a^{(j)}_i$
or $a^{(j)}_i+1$, $a^{(j)}_i<a^{(j)}_{i+1}$, $a^{(1)}_i=i$ or $i+1\}$.

Step~3.
One-to-one correspondence between paths and terms of $\Delta^L(k;{\bf i})(\tau)$.

Let us recall that $(0;a^{(m)}_1,a^{(m)}_2,\dots,a^{(m)}_d)$ means a~minor of the identity matrix, and\linebreak
$(0;a^{(m)}_1,a^{(m)}_2,\dots,a^{(m)}_d)=\det(D(0;a^{(m)}_1),D(0;a^{(m)}_2),\dots,D(0;a^{(m)}_d))$, and each $D(0;i)$ is
given as follows:
\begin{gather*}
D(0;i)=
\begin{pmatrix}
x^{(0)}_{m'+1,i}
\\
\vdots
\\
x^{(0)}_{m'+d,i}
\end{pmatrix}
=
\begin{pmatrix}
\delta_{m'+1,i}
\\
\vdots
\\
\delta_{m'+d,i}
\end{pmatrix}
.
\end{gather*}
By $1\leq a^{(m)}_1<a^{(m)}_2<\dots<a^{(m)}_d \leq m+d$ (Def\/inition~\ref{pathdef}), we have
\begin{gather*}
\big(0;a^{(m)}_1,a^{(m)}_2,\dots,a^{(m)}_d\big)=
\begin{cases}
1 & \text{if} \quad  a^{(m)}_1=m'+1,\  a^{(m)}_2=m'+2,\ \dots,\  a^{(m)}_d=m'+d,
\\
0 & \text{otherwise}.
\end{cases}
\end{gather*}
Therefore, $\Delta^L(k;{\bf i})(\tau)=\det(D(m;1),D(m;2),\dots,D(m;d))$ is equal to the coef\/f\/icient of
$(0;m'+1$, $m'+2,\dots,m'+d)$.
By~\eqref{lincom1}, setting $a^{(0)}_i=i$, $a^{(m)}_i=m'+i$, we have
\begin{gather}
\Delta^L(k;{\bf i})(\tau)=\sum\limits_{a^{(j)}_i,\ i\in \{1,\dots,d\},\ 1\leq j\leq m-1}
\left(\prod\limits^d_{i=1}\frac{\tau_{l_{m-1}+a^{(1)}_i-1}}{\tau_{l_{m-1} +a^{(0)}_i}}\right.
\nonumber
\\
\phantom{\Delta^L(k;{\bf i})(\tau)=}
\times\left.\prod\limits^d_{i=1}
\frac{\tau_{l_{m-2}+a^{(2)}_i-1}}{\tau_{l_{m-2}+a^{(1)}_i}} \times\dots\times\prod\limits^d_{i=1}\frac{\tau_{l_{0}
+a^{(m)}_i-1}}{\tau_{l_{0}+a^{(m-1)}_i}} \right),
\label{summand1}
\end{gather}
where the index $a^{(j)}_i$ runs over $\{(a^{(j)}_i:i\in \{1,\dots,d\}$, $1\leq j\leq m-1) \,|\, a^{(j+1)}_i=a^{(j)}_i$
or $a^{(j)}_i+1$, $a^{(j)}_i<a^{(j)}_{i+1}$, $0\leq j\leq m-1\}$.
Note that these conditions equal to the ones in Def\/inition~\ref{pathdef}(ii) and (iii), and the conditions
$a^{(0)}_i=i$, $a^{(m)}_i=m'+i$ are equal to the ones in Def\/inition~\ref{pathdef}(iv).
We see that the summand in~\eqref{summand1} is equal to the label of the path
\begin{gather*}
\big(m;a^{(0)}_1,\dots,a^{(0)}_d\big)\rightarrow\big(m-1;a^{(1)}_1,\dots,a^{(1)}_d\big)\rightarrow\big(m-2;a^{(2)}_1,\dots,a^{(2)}_d\big)\rightarrow
\cdots
\\
\hphantom{\big(m;a^{(0)}_1,\dots,a^{(0)}_d\big)}{}
\rightarrow\big(1;a^{(m-1)}_1,\dots,a^{(m-1)}_d\big)\rightarrow\big(0;a^{(m)}_1,\dots,a^{(m)}_d\big)\in X_d(m,m').
\end{gather*}
Hence we obtain $\Delta^L(k;{\bf i})(\tau)=\sum\limits_{p\in X_d(m,m')} Q(p)$.
\end{proof}

\begin{ex}
\label{pathex2}
We use the same notations~$m$, $m'$ and $d=i_k$ at the beginning of Section~\ref{gmc}.
We set rank $r=4$, $u=s_1s_2s_3s_4s_1s_2s_3s_1s_2s_1\in W$ and $k=6$, which is the same setting as in introduction.
Let
\begin{gather*}
{\bf i}:=(1,2,3,4,1,2,3,1,2,1)
\end{gather*}
be a~reduced for~$u$, and let
\begin{gather*}
{\bf i}'=(\underbrace{1,2,3,4}_{\text{$1^{\rm st}$ cycle}},\underbrace{1,2,3}_{\text{$2^{\rm nd}$ cycle}},\underbrace{1,2}_{\text{$3^{\rm rd}$ cycle}})
\end{gather*}
be a~reduced word for $u s_1$.
By Lemma~\ref{gmlem}, we get $\Delta^L(6;{\bf i})(\tau)=\Delta^L(6;{\bf i}')(\tau)$.
Since $m=3$, $m'=2$ and $d:=i_k=2$, it follows from Example~\ref{pathex} and Proposition~\ref{pathlem} that
\begin{gather*}
\Delta^L(6;{\bf i})(\tau) = \Delta^L(6;{\bf i}')(\tau)=\sum\limits_{p\in X_2(3,2)} Q(p)
\\
\phantom{\Delta^L(6;{\bf i})(\tau)}
 =  \frac{1}{\tau_{l_{2}+2}}+\frac{\tau_{l_1+2}}{\tau_{l_2+1}\tau_{l_1+3}}
+\frac{\tau_{l_0+3}}{\tau_{l_2+1}\tau_{l_0+4}}+\frac{\tau_{l_1+1}}{\tau_{l_1+3}}
+\frac{\tau_{l_1+1}\tau_{l_0+3}}{\tau_{l_1+2}\tau_{l_0+4}} +\frac{\tau_{l_0+2}}{\tau_{l_0+4}}
\\
\phantom{\Delta^L(6;{\bf i})(\tau)}
 =  \frac{1}{\tau_{9}}+\frac{\tau_{6}}{\tau_{8}\tau_{7}} +\frac{\tau_{3}}{\tau_{8}\tau_{4}}+\frac{\tau_{5}}{\tau_{7}}
+\frac{\tau_{5}\tau_{3}}{\tau_{6}\tau_{4}} +\frac{\tau_{2}}{\tau_{4}},
\end{gather*}
which is equal to \eqref{Del-fin}.
\end{ex}

\subsection[The explicit description of $\Delta^L(k;{\bf i})(\tau)$]{The explicit description of $\boldsymbol{\Delta^L(k;{\bf i})(\tau)}$}

In Proposition~\ref{pathlem}, we had described the terms of $\Delta^L(k;{\bf i})(\tau)$ as the paths in $X_d(m,m')$.
In this subsection, we shall describe $\Delta^L(k;{\bf i})(\tau)$ explicitly by using some properties of paths.
This description will be used in the proof of Theorem~\ref{thm2}.
First, we need to show some lemmas.

Let us write a~path $p\in X_d(m,m')$ as follows
\begin{gather}
p=\big(m;a^{(0)}_1,\dots,a^{(0)}_d\big)\rightarrow\big(m-1;a^{(1)}_1,\dots,a^{(1)}_d\big)\rightarrow
\big(m-2;a^{(2)}_1,\dots,a^{(2)}_d\big)\rightarrow\cdots
\nonumber\\
\hphantom{p=}{}
\rightarrow\big(1;a^{(m-1)}_1,\dots,a^{(m-1)}_d\big)\rightarrow\big(0;a^{(m)}_1,\dots,a^{(m)}_d\big).
\label{fixpath}
\end{gather}

We begin with the following lemma:
\begin{lem}
\label{nmblem}
For a~path~$p$ \eqref{fixpath} and $i\in \{1,\dots,d\}$, we have
\begin{gather*}
\#\big\{s \,|\, a^{(s)}_i=a^{(s+1)}_i,\  0\leq s\leq m-1 \big\}=m-m'.
\end{gather*}
\end{lem}

\begin{proof} By Def\/inition~\ref{pathdef}(iii) and (iv), we have
\begin{gather*}
i=a^{(0)}_i\leq a^{(1)}_i\leq\dots\leq a^{(m)}_i=m'+i,
\qquad
a^{(s+1)}_i=a^{(s)}_i
\qquad
\text{or}
\qquad
a^{(s)}_i+1.
\end{gather*}
Thus, we get
\begin{gather*}
\#\big\{s \,|\, a^{(s+1)}_i=a^{(s)}_i+1,\ 0\leq s\leq m-1\big\}=m',
\end{gather*}
which implies that $\#\{s \,|\, a^{(s)}_i=a^{(s+1)}_i,\ 0\leq s\leq m-1 \}=m-m'$.
\end{proof}

\begin{defn}
\label{sijdefi}
For a~path~$p$ \eqref{fixpath} and $i\in \{1,\dots,d\}$, we set $\{q^{(j)}_i\}_{1\leq j\leq m-m'}$,
$0\leq q^{(1)}_i<\dots<q^{(m-m')}_i\leq m-1$, as
\begin{gather}
\label{ldef}
\big\{q^{(1)}_i,\; q^{(2)}_i,\dots,q^{(m-m')}_i\big\}:=\big\{q \,|\, a^{(q)}_i=a^{(q+1)}_i,\;
0\leq q\leq m-1 \big\}.
\end{gather}
We also set $k^{(j)}_i\in [1,m'+d]$ ($1\leq i\leq d$, $1\leq j\leq m-m'$) as
\begin{gather}
\label{kdef}
k^{(j)}_i:=a^{(q^{(j)}_i)}_i.
\end{gather}
\end{defn}

\begin{lem}\label{kpro2}\quad
\begin{enumerate}\itemsep=0pt
\item[$(i)$] For $1\leq i\leq d$ and $1\leq j\leq m-m'$,
\begin{gather*}
q^{(j)}_{i}=k^{(j)}_{i}+j-i-1.
\end{gather*}
\item[$(ii)$] For $1\leq j\leq m-m'$ and $1\leq i\leq d-1$,
\begin{gather*}
1\leq k^{(j)}_{i}<k^{(j)}_{i+1}\leq m'+d,
\qquad
q^{(j)}_i\leq q^{(j)}_{i+1}.
\end{gather*}
\end{enumerate}
\end{lem}

\begin{proof}
(i) The def\/inition of $q^{(j)}_i$ in~\eqref{ldef} means that the path~$p$ has the following~$i$-sequence
(Def\/inition~\ref{iseq}):
\begin{gather}
a^{(0)}_i=i, \quad a^{(1)}_i=i+1, \quad a^{(2)}_i=i+2, \quad \dots, \quad a^{(q^{(1)}_i)}_i=i+q^{(1)}_i,\nonumber\\
a^{(q^{(1)}_i+1)}_i=i+q^{(1)}_i, \quad a^{(q^{(1)}_i+2)}_i=i+q^{(1)}_i+1, \quad \dots, \quad a^{(q^{(2)}_i)}_i=i+q^{(2)}_i-1,\nonumber\\
a^{(q^{(2)}_i+1)}_i=i+q^{(2)}_i-1, \quad a^{(q^{(2)}_i+2)}_i=i+q^{(2)}_i,\quad \dots, \quad a^{(q^{(3)}_i)}_i=i+q^{(3)}_i-2,\nonumber\\
\label{jlist}
\cdots\cdots\cdots\cdots\cdots\cdots\cdots\cdots\cdots\cdots\cdots\cdots\cdots\cdots\cdots\cdots\cdots\cdots\cdots\cdots\cdots\cdots\cdots\cdots
\\
\nonumber
a^{(q^{(j-1)}_i+1)}_i=i+q^{(j-1)}_i-j+2,
\quad
a^{(q^{(j-1)}_i+2)}_i=i+q^{(j-1)}_i-j+3,\quad
\dots,\\
a^{(q^{(j)}_i)}_i=i+q^{(j)}_i-j+1,\quad
a^{(q^{(j)}_i+1)}_i=i+q^{(j)}_i-j+1, \quad a^{(q^{(j)}_i+2)}_i=i+q^{(j)}_i-j+2, \quad \dots.\nonumber
\end{gather}

Hence we have
\begin{gather*}%\label{kpro2pr1}
k^{(j)}_{i}=a^{(q^{(j)}_{i})}_i=i+q^{(j)}_{i}-j+1,
\end{gather*}
which implies $q^{(j)}_{i}=k^{(j)}_{i}+j-i-1$.

(ii) By Def\/inition~\ref{pathdef}(iii),
\begin{gather}
\label{kpro2pr3ano}
i=a^{(0)}_{i}\leq a^{(1)}_{i}\leq\dots\leq a^{(q^{(j)}_{i+1})}_{i}\leq a^{(q^{(j)}_{i+1}+1)}_{i},
\\
\nonumber
a^{(\zeta)}_{i}=a^{(\zeta-1)}_{i}
\qquad
\text{or}
\qquad
a^{(\zeta-1)}_{i}+1,
\qquad
1\leq\zeta\leq q^{(j)}_{i+1}+1.
\end{gather}
We obtain
\begin{gather}
\label{kpro2pr4ano}
q^{(j)}_{i+1}+1-j\geq\#\big\{\zeta \,|\,  a^{(\zeta)}_{i}=a^{(\zeta-1)}_{i}+1,\ 1\leq\zeta\leq q^{(j)}_{i+1}+1\big\},
\end{gather}
otherwise, it follows from~\eqref{kpro2pr3ano} and (i)
that $a^{(q^{(j)}_{i+1}+1)}_i>i+q^{(j)}_{i+1}+1-j=k_{i+1}-1=a^{(q^{(j)}_{i+1})}_{i+1}-1$,
and hence $a^{(q^{(j)}_{i+1}+1)}_i\geq a^{(q^{(j)}_{i+1})}_{i+1}=a^{(q^{(j)}_{i+1}+1)}_{i+1}$, which contradicts Def\/inition~\ref{pathdef}(ii).

The inequality~\eqref{kpro2pr4ano} means that
\begin{gather}
\label{kpro2pr5ano}
j\leq\#\big\{\zeta \,|\,  a^{(\zeta)}_{i}=a^{(\zeta-1)}_{i},\ 1\leq\zeta\leq q^{(j)}_{i+1}+1\big\}.
\end{gather}
On the other hand, using the list~\eqref{jlist} (or the def\/inition of $q^{(j)}_i$ in~\eqref{ldef}), we have
\begin{gather}
\label{kpro2pr6ano}
j=\#\big\{\zeta \,|\,  a^{(\zeta)}_{i}=a^{(\zeta-1)}_{i},\ 1\leq\zeta\leq q^{(j)}_{i}+1\big\}.
\end{gather}
Since $a^{(q^{(j)}_i)}_i=a^{(q^{(j)}_i+1)}_i$, the equation~\eqref{kpro2pr6ano} means
\begin{gather}
\label{kpro2pr6ano2}
j-1=\#\big\{\zeta \,|\,  a^{(\zeta)}_{i}=a^{(\zeta-1)}_{i},\ 1\leq\zeta\leq q^{(j)}_{i}\big\}.
\end{gather}
Thus, by~\eqref{kpro2pr5ano} and~\eqref{kpro2pr6ano2}, we have $q^{(j)}_{i}< q^{(j)}_{i+1}+1$, and hence
$q^{(j)}_{i}\leq q^{(j)}_{i+1}$, which yields $k^{(j)}_i<k^{(j)}_{i+1}$ since
$k^{(j)}_i=i+q^{(j)}_i-j+1<i+q^{(j)}_{i+1}-j+2=(i+1)+q^{(j)}_{i+1}-j+1=k^{(j)}_{i+1}$.
Remark~\ref{pathrem} implies that $k^{(j)}_{i+1}:=a^{(q^{(j)}_{i+1})}_{i+1}\leq m'+d$.
\end{proof}

For $1\leq i\leq m$ and $1\leq j\leq r$, we set the Laurent monomials
\begin{gather}
\label{ccbar}
\overline{C}(i,j):=\frac{\tau_{l_i+j-1}}{\tau_{l_i+j}}.
\end{gather}

\begin{lem}
\label{thm1lem}
For a~path~$p$ \eqref{fixpath}, we set $q^{(j)}_i$ and $k^{(j)}_i$ as in~\eqref{ldef} and~\eqref{kdef}.
Then we have
\begin{gather}
\label{thm1lemclaim}
Q(p)=\prod\limits^{d}_{i=1} \prod\limits_{j=1}^{m-m'}\overline{C}\big(m-q^{(j)}_i-1,k^{(j)}_i\big).
\end{gather}
\end{lem}

\begin{proof}
Let us recall the def\/inition of $Q(p)$ (Def\/inition~\ref{labeldef}(ii)):
\begin{gather}
\label{qpprod}
Q(p)=\prod\limits_{s=0}^{m-1}\left(\prod\limits^d_{i=1}\frac{\tau_{l_{m-s-1}+a^{(s+1)}_i-1}}{\tau_{l_{m-s-1}+a^{(s)}_i}}
\right).
\end{gather}

For each $i=1,2,\dots,d$, the~$i$-sequence
\begin{gather*}
i=a^{(0)}_i\leq a^{(1)}_i\leq a^{(2)}_i\leq \dots \leq a^{(m)}_i=m'+i
\end{gather*}
of~$p$ satisf\/ies $a^{(s+1)}_i=a^{(s)}_i$ or $a^{(s)}_i+1$, $0\leq s\leq m-1$, by Def\/inition~\ref{pathdef}(iii).
If $a^{(s+1)}_i=a^{(s)}_i+1$, then
\begin{gather*}
\frac{\tau_{l_{m-s-1}+a^{(s+1)}_i-1}}{\tau_{l_{m-s-1}+a^{(s)}_i}}=1.
\end{gather*}
Therefore, it follows from the def\/inition~\eqref{ldef} of $q^{(j)}_i$ and
$k^{(j)}_i:=a^{(q^{(j)}_i)}_i=a^{(q^{(j)}_i+1)}_i$ that
\begin{gather*}
\prod\limits_{s=0}^{m-1}\frac{\tau_{l_{m-s-1}+a^{(s+1)}_i-1}}{\tau_{l_{m-s-1}+a^{(s)}_i}}  =
\prod\limits_{j=1}^{m-m'}\frac{\tau_{l_{m-q^{(j)}_i-1}+k^{(j)}_i-1}}{\tau_{l_{m-q^{(j)}_i-1}+k^{(j)}_i}}
 = \prod\limits_{j=1}^{m-m'}\overline{C}\big(m-q^{(j)}_i-1,k^{(j)}_i\big),
\end{gather*}
which implies~\eqref{thm1lemclaim} by~\eqref{qpprod}.
\end{proof}

Let us describe $\Delta^L(k;{\bf i})(\tau)$ explicitly.

\begin{prop}
\label{monoprop}
\begin{gather*}
\Delta^L(k;{\bf i})(\tau)= \sum\limits_{(*)}\prod\limits^{d}_{i=1}
\prod\limits_{j=1}^{m-m'}\overline{C}\big(m-K^{(j)}_i-j+i,K^{(j)}_i\big).
\end{gather*}
where $(*)$ is the conditions for $K^{(j)}_i$ $(1\leq i\leq d$, $1\leq j\leq m-m')$:
$1\leq K^{(j)}_1<K^{(j)}_2<\dots<K^{(j)}_d\leq m'+d$, $1\leq j\leq m-m'$,
$i\leq K^{(1)}_i\leq K^{(2)}_i\leq\dots\leq K^{(m-m')}_i\leq m'+i$, $1\leq i\leq d$.
\end{prop}

\begin{proof}
Using Lemmas~\ref{kpro2}(i) and~\ref{thm1lem}, we see that $Q(p)$, $p\in X_d(m,m')$, is described~as
\begin{gather*}
Q(p)=\prod\limits^{d}_{i=1} \prod\limits_{j=1}^{m-m'}\overline{C}\big(m-k^{(j)}_i-j+i,k^{(j)}_i\big).
\end{gather*}
with $\{k^{(j)}_i\}_{1\leq i\leq d, 1\leq j\leq m-m'}$ which satisfy the conditions in Lemma~\ref{kpro2}(ii), that is,
$1\leq k^{(j)}_1<k^{(j)}_2<\dots<k^{(j)}_d\leq m'+d$.
Furthermore, Def\/inition~\ref{pathdef}(iii) and~(iv) show that $i=a^{(0)}_i\leq a^{(1)}_i\leq \cdots\leq a^{(m)}_i=m'+i$,
which means that $i\leq k^{(1)}_i\leq k^{(2)}_i\leq\cdots\leq k^{(m-m')}_i\leq m'+i$ for $1\leq i\leq d$~by
$k^{(j)}_i:=a^{(q^{(j)}_i)}_i$ and $q^{(j)}_i<q^{(j+1)}_i$ (see Def\/inition~\ref{sijdefi}).
Thus, $\{k^{(j)}_i\}$ satisf\/ies the conditions $(*)$ in Proposition~\ref{monoprop}.

Conversely, let $\{K^{(j)}_i\}_{1\leq i\leq d, 1\leq j\leq m-m'}$ the set of numbers which satisf\/ies the conditions
$(*)$ in Proposition~\ref{monoprop}:
\begin{gather*}
%\label{cond1}
1\leq K^{(j)}_1<K^{(j)}_2<\dots<K^{(j)}_d\leq m'+d,
\qquad
1\leq j\leq m-m',
\end{gather*}
and
\begin{gather}
\label{cond2}
i\leq K^{(1)}_i\leq \dots\leq K^{(m-m')}_i\leq m'+i,
\qquad
1\leq i\leq d.
\end{gather}
We set
\begin{gather}
\label{QDEF}
Q^{(j)}_{i}:= K^{(j)}_{i}+j-i-1,
\end{gather}
for $1\leq i\leq d$ and $1\leq j\leq m-m'$.
We need to show that there exists a~path $p\in X_d(m,m')$ such that
\begin{gather}
\label{finclaim}
Q(p)=\prod\limits^{d}_{i=1} \prod\limits_{j=1}^{m-m'}\overline{C}\big(m-Q^{(j)}_i-1,K^{(j)}_i\big)=\prod\limits^{d}_{i=1}
\prod\limits_{j=1}^{m-m'}\overline{C}\big(m-K^{(j)}_i-j+i,K^{(j)}_i\big).
\end{gather}

Since we supposed $K^{(j)}_i<K^{(j)}_{i+1}$, we can easily verify
\begin{gather}
\label{lequ1}
Q^{(j)}_i\leq Q^{(j)}_{i+1},
\end{gather}
by~\eqref{QDEF}.
We claim that $0\leq Q^{(j)}_i\leq m-1$ for $1\leq i\leq d$ and $1\leq j\leq m-m'$.
By the condition~\eqref{cond2}, we get $i\leq K^{(j)}_i$.
So it is clear that $0\leq Q^{(j)}_i$.
It follows from the condition~\eqref{cond2} and~\eqref{QDEF} that $Q^{(j)}_{i}=K^{(j)}_{i}+j-i-1\leq
m'+i+j-i-1=m'+j-1\leq m-1$.
Therefore, we have $0\leq Q^{(j)}_i\leq m-1$ for all $1\leq i\leq d$ and $1\leq j\leq m-m'$.

We def\/ine a~path $p=(m;a^{(0)}_1,\dots,a^{(0)}_d)\rightarrow\dots \rightarrow (0;a^{(m)}_1,\dots,a^{(m)}_d)\in
X_d(m,m')$ as follows. For~$i$, $1\leq i\leq d$, we def\/ine the~$i$-sequence (Def\/inition~\ref{iseq}) of~$p$ as
\begin{gather}
a^{(0)}_i=i,
\quad
a^{(1)}_i=i+1,
\quad
a^{(2)}_i=i+2,\quad \dots, \quad a^{(Q^{(1)}_i)}_i=i+Q^{(1)}_i,
\nonumber
\\
a^{(Q^{(1)}_i+1)}_i=i+Q^{(1)}_i,
\quad
a^{(Q^{(1)}_i+2)}_i=i+Q^{(1)}_i+1,\quad \dots,\quad a^{(Q^{(2)}_i)}_i=i+Q^{(2)}_i-1,
\nonumber
\\
a^{(Q^{(2)}_i+1)}_i=i+Q^{(2)}_i-1,
\quad
a^{(Q^{(2)}_i+2)}_i=i+Q^{(2)}_i,\quad \dots,\quad a^{(Q^{(3)}_i)}_i=i+Q^{(3)}_i-2,
\nonumber
\\
\label{jlist2}
\cdots\cdots\cdots\cdots\cdots\cdots\cdots\cdots\cdots\cdots\cdots\cdots\cdots\cdots\cdots\cdots\cdots\cdots
\\
a^{(Q^{(m-m'-1)}_i+1)}_i=i+Q^{(m-m'-1)}_i-m+m'+2,\quad \dots,
\nonumber
\\
a^{(Q^{(m-m')}_i)}_i=i+Q^{(m-m')}_i-m+m'+1,
\qquad
a^{(Q^{(m-m')}_i+1)}_i=i+Q^{(m-m')}_i-m+m'+1,
\nonumber
\\
a^{(Q^{(m-m')}_i+2)}_i=i+Q^{(m-m')}_i-m+m'+2,
\nonumber
\\
a^{(Q^{(m-m')}_i+3)}_i=i+Q^{(m-m')}_i-m+m'+3,\quad \dots, \quad a^{(m)}_i=m'+i.
\nonumber
\end{gather}

It is easy to see that $a^{(Q^{(j)}_i)}_i=K^{(j)}_i$, $1\leq j\leq m-m'$, by~\eqref{QDEF} and~\eqref{jlist2}.
Clearly, the path~$p$ satisf\/ies Def\/inition~\ref{pathdef}(iii) and~(iv).
For $1\leq s\leq m$, we obtain $a^{(s)}_i<a^{(s)}_{i+1}$ by~\eqref{lequ1} and~\eqref{jlist2}.

Hence~$p$ is well-def\/ined, and~\eqref{finclaim} is follows from Lemma~\ref{thm1lem} (see~\eqref{ldef} and~\eqref{kdef}).
Thus, Proposition~\ref{monoprop} follows from Proposition~\ref{pathlem}.
\end{proof}

\begin{ex}
\label{pathex3}
We use the same setting in Example~\ref{pathex2}.
Since $m=3$, $m'=2$ and $d:=i_k=2$, it follows from Proposition~\ref{monoprop} that
\begin{gather*}
\Delta^L(6;{\bf i})(\tau) =  \sum\limits_{1\leq K_1<K_2\leq 4} \left(\prod\limits^{2}_{i=1}
\overline{C}(i-K_{i}+2,K_{i})\right)
\\
\phantom{\Delta^L(6;{\bf i})(\tau)}
 =  \overline{C}(2,1)\overline{C}(2,2)+\overline{C}(2,1)\overline{C}(1,3)+\overline{C}(2,1)\overline{C}(0,4)
\\
\phantom{\Delta^L(6;{\bf i})(\tau) =}
{} +\overline{C}(1,2)\overline{C}(1,3)+\overline{C}(1,2)\overline{C}(0,4)+\overline{C}(0,3)\overline{C}(0,4)
\\
\phantom{\Delta^L(6;{\bf i})(\tau)}
 =  \frac{1}{\tau_{l_{2}+2}}+\frac{\tau_{l_1+2}}{\tau_{l_2+1}\tau_{l_1+3}}
+\frac{\tau_{l_0+3}}{\tau_{l_2+1}\tau_{l_0+4}}+\frac{\tau_{l_1+1}}{\tau_{l_1+3}}
+\frac{\tau_{l_1+1}\tau_{l_0+3}}{\tau_{l_1+2}\tau_{l_0+4}} +\frac{\tau_{l_0+2}}{\tau_{l_0+4}},
\end{gather*}
which is equal to the one in Example~\ref{pathex2}.
\end{ex}

\subsection{The completion of the proof of Theorem~\ref{thm2}}\label{prthm2fin}

 In this subsection, we shall complete the proof of Theorem~\ref{thm2}.

Let us recall the def\/inition~\eqref{asidef} of $A_{s,i}$.
Since we identify the variables $\{Y_{s,i}\}$ with $\{\tau_{l_s+i}\}$, we have the following (see Remark~\ref{tautau})
\begin{gather*}
A_{s,i}=\frac{\tau_{l_s+i}\tau_{l_{s+1}+i}}{\tau_{l_s+i+1}\tau_{l_{s+1}+i-1}},
\qquad
0\leq s\leq m-1,
\qquad
1\leq i\leq r.
\end{gather*}
Therefore, by the def\/inition~\eqref{ccbar} of $\overline{C}(s,i)$, we have
\begin{gather}
\label{canda}
\overline{C}(s,i)\cdot A_{s-1,i}=\overline{C}(s-1,i+1).
\end{gather}

\begin{lem}
\label{pathlems}
For a~path $p\in X_d(m,m')$, we describe the monomial $Q(p)$ as in Lemma~{\rm \ref{thm1lem}}
\begin{gather}
\label{qpexpand}
Q(p)=\prod\limits^{d}_{i=1} \prod\limits_{j=1}^{m-m'}\overline{C}\big(m-K^{(j)}_i-j+i,K^{(j)}_i\big),
\end{gather}
where $\{K^{(j)}_i\}$ satisfies the condition $(*)$ in Proposition~{\rm \ref{monoprop}}.
For~$s$, $1\leq s<m'+d$, if $\tilde{e}_s Q(p)\neq 0$, then there exist~$i$ and~$j$, $1\leq i\leq d$, $1\leq j\leq m-m'$,
such that
\begin{gather}
\label{pathlems1}
K^{(j)}_i=s
\qquad
\text{and}
\qquad
\tilde{e}_s Q(p)=Q(p)\cdot A_{m-K^{(j)}_i-j+i-1,s}.
\end{gather}
Furthermore, there exists a~path $p'\in X_d(m,m')$ such that
\begin{gather}
\label{pathlems2}
\tilde{e}_s Q(p)=Q(p').
\end{gather}
\end{lem}

\begin{proof}
We suppose that the monomial $Q(p)$ does not include factor as in the form $\tau_{l_t+s}^{-1}$, $0\leq t\leq m-1$,
which means wt$(Q(p))(h_i)=\varphi_i(Q(p))$.
Hence $\varepsilon_i(Q(p)):=\varphi_i(Q(p))-\text{wt}(Q(p))(h_i)=0$, which contradicts the assumption $\tilde{e}_s
Q(p)\neq 0$.
Hence the monomial $Q(p)$ includes factor $\tau_{l_t+s}^{-1}$, $0\leq t\leq m-1$.
We set the numbers $0\leq t_1<t_2<\dots<t_{\xi}\leq m-1$~by
\begin{gather*}
\{t_1,t_2,\dots,t_{\xi}\}:=\big\{t \,|\, Q(p)\ \text{includes\ factors}\ \tau_{l_t+s}^{-1}\big\},
\qquad
1\leq\xi.
\end{gather*}
The def\/inition~\eqref{ccbar} of $\overline{C}(t,s)$ and~\eqref{qpexpand} show that there exist $1\leq
i_1,\dots,i_{\xi}\leq d$ and $1\leq j_1,\dots,j_{\xi}\leq m-m'$ such that
\begin{gather}
\label{ktqeq}
K^{(j_a)}_{i_a}=s
\qquad
\text{and}
\qquad
t_a=m-K^{(j_a)}_{i_a}-j_a+i_a,
\qquad
a=1,2,\dots,\xi.
\end{gather}

As in Example~\ref{exex}, for a~given $s\in I$ and monomial $Y=\prod\limits_{q \in \mathbb{Z},\  i \in I}
\tau_{l_q+i}^{\zeta_{q,i}}$, we def\/ine $\nu_Y(n):=\sum\limits_{q\leq n}\zeta_{q,s}$, $n\in \mathbb{Z}$.

We set $Y:=Q(p)$.
We claim that $n_{e_s}\in \{t_1-1,\dots,t_{\xi}-1\}$~\eqref{nxi}, otherwise, we get $n_{e_s}\notin
\{t_1-1,\dots,t_{\xi}-1\}$.
Then $\varphi_s(Y)=\nu_Y(n_{e_s})\leq \nu_Y(n_{e_s}+1)$.
If $\nu_Y(n_{e_s})< \nu_Y(n_{e_s}+1)$, then it contradicts the def\/inition of
\begin{gather}
\label{nesta0}
\varphi_s(Y):=\max \{\nu_Y(n) \,|\, n\in\mathbb{Z}\}.
\end{gather}
If $\nu_Y(n_{e_s})=\nu_Y(n_{e_s}+1)$, then it contradicts the def\/inition of $n_{e_s}:=\max \{n \,|\, \varphi_s(Y)=\nu_Y(n)
\}$.
Hence we obtain $n_{e_s}\in \{t_1-1,\dots,t_{\xi}-1\}$, which implies there exists~$a$, $1\leq a\leq \xi$, such that
$n_{e_s}=t_a-1$.
By~\eqref{ktqeq}, we have $n_{e_s}=t_a-1=m-K^{(j_a)}_{i_a}-j_a+i_a-1$.
Therefore, we obtain $\tilde{e}_s Y:= Y\cdot A_{n_{e_s},s}=Y\cdot A_{m-K^{(j_a)}_{i_a}-j_a+i_a-1,s}$.

Since $Q(p)$ includes the factor $\tau_{l_{t_a}+s}^{-1}$, the path~$p$ includes the following edge
\begin{gather*}
%\begin{xy}
%(0,25) *{(t_a+1;a^{(m-t_a-1)}_1,\cdots,a^{(m-t_a-1)}_{q-1},\underset{q\ th}{s},a^{(m-t_a-1)}_{q+1},\cdots,a^{(m-t_a-1)}_{d})}="(s+1,a)",
%(0,12)*{(t_a;a^{(m-t_a)}_1,\cdots,a^{(m-t_a)}_{q-1},\underset{q\ th}{s},a^{(m-t_a)}_{q+1},\cdots,a^{(m-t_a)}_{d})}="(s,a+1)",
%(0,0)*{(t_a-1;a^{(m-t_a+1)}_1,\cdots,a^{(m-t_a+1)}_{q-1},a^{(m-t_a+1)}_{q},a^{(m-t_a+1)}_{q+1},\cdots,a^{(m-t_a+1)}_{d})}="(s-1,a+1)",
%\ar@{->} "(s+1,a)";"(s,a+1)"
%\ar@{->} "(s,a+1)";"(s-1,a+1)"
%\end{xy},
%\begin{xy}
\big(t_a+1;a^{(m-t_a-1)}_1,\cdots,a^{(m-t_a-1)}_{q-1},\underset{q^{\rm th}}{s},a^{(m-t_a-1)}_{q+1},\cdots,a^{(m-t_a-1)}_{d}\big)\\
\qquad{}\rightarrow \big(t_a;a^{(m-t_a)}_1,\cdots,a^{(m-t_a)}_{q-1},\underset{q^{\rm th}}{s},a^{(m-t_a)}_{q+1},\cdots,a^{(m-t_a)}_{d}\big)\\
\qquad{}\rightarrow
\big(t_a-1;a^{(m-t_a+1)}_1,\cdots,a^{(m-t_a+1)}_{q-1},a^{(m-t_a+1)}_{q},a^{(m-t_a+1)}_{q+1},\cdots,a^{(m-t_a+1)}_{d}\big)
\end{gather*}
for some $1\leq q\leq d$.
We claim that $a^{(m-t_a)}_{q+1}>s+1$, otherwise $a^{(m-t_a)}_{q+1}=s+1$ and $a^{(m-t_a-1)}_{q+1}=s+1$~by
Def\/inition~\ref{pathdef}(ii) and~(iii), which implies that $Q(p)$ does not include the factor $\tau_{l_{t_a}+s}^{-1}$
by Def\/inition~\ref{labeldef}(i).
It contradicts the def\/inition of $t_a$.
So we get
\begin{gather}
\label{lem1pr0}
a^{(m-t_a)}_{q+1}>s+1.
\end{gather}
If $t_a=0$, then Def\/inition~\ref{pathdef}(iv) means that $s=m'+q$ (by the assumption $s<m'+d$, we get $q<d$), and
$a^{(m-t_a)}_{q+1}=m'+q+1=s+1$, which contradicts~\eqref{lem1pr0}.
So we have $1\leq t_a$.
We also claim that $a^{(m-t_a+1)}_{q}>s$, otherwise $a^{(m-t_a+1)}_{q}=s$, which implies that $Q(p)$ includes the factor
$\tau_{l_{t_a-1}+s}^{-1}$ by $a^{(m-t_a+1)}_{q+1}\geq a^{(m-t_a)}_{q+1}>s+1$.
Therefore, $\varphi_s(Y)=\nu_Y(n_{e_s})=\nu_Y(t_a-1)<\nu_Y(t_a-2)$, which contradicts~\eqref{nesta0}.
So we get $a^{(m-t_a+1)}_{q}>s$, which means
\begin{gather}
\label{lem1pr1}
a^{(m-t_a+1)}_{q}=s+1,
\end{gather}
by Def\/inition~\ref{pathdef}(iii).

Let $p'\in X_d(m,m')$ be the path obtained from~$p$ by replacing the vertex
\begin{gather*}
\big(t_a;a^{(m-t_a)}_1,\dots,a^{(m-t_a)}_{q-1},\underset{q^{\rm th}}{s},a^{(m-t_a)}_{q+1},\dots,a^{(m-t_a)}_{d}\big)
\end{gather*}
by
\begin{gather*}
\big(t_a;a^{(m-t_a)}_1,\dots,a^{(m-t_a)}_{q-1},\underset{q^{\rm th}}{s+1},a^{(m-t_a)}_{q+1},\dots,a^{(m-t_a)}_{d}\big).
\end{gather*}
By~\eqref{lem1pr0} and~\eqref{lem1pr1}, the path $p'$ is well-def\/ined.
Def\/inition~\ref{labeldef}(i) shows that
\begin{gather*}
Q(p')=Q(p)\cdot \frac{\tau_{l_{t_a-1}+s}\tau_{l_{t_a}+s}}{\tau_{l_{t_a-1}+s+1}\tau_{l_{t_a}+s-1}}=Q(p)\cdot
A_{m-K^{(j_a)}_{i_a}-j_a+i_a-1,s}=\tilde{e}_s Q(p)
\end{gather*}
by~\eqref{ktqeq}.
\end{proof}

Next lemma shows that the coef\/f\/icients $c(b)$ in Theorem~\ref{thm2} are equal to $1$ for all $b\in B^-_{u_{\leq
k}}((m'-m)\Lambda_d)$ by Proposition~\ref{pathlem}.

\begin{lem}
\label{coefflem}
For paths $p,p'\in X_d(m,m')$, if $p\neq p'$ then $Q(p)\neq Q(p')$.
\end{lem}

\begin{proof}
We suppose that $Q(p)=Q(p')$.
Let us prove $p=p'$.
We denote~$p$ and $p'$~by
\begin{gather*}
p=\big(m;a^{(0)}_1,\dots,a^{(0)}_d\big)\rightarrow\cdots\rightarrow \big(0;a^{(m)}_1,\dots,a^{(m)}_d\big)
\end{gather*}
and
\begin{gather*}
p'=\big(m;b^{(0)}_1,\dots,b^{(0)}_d\big)\rightarrow\dots\rightarrow \big(0;b^{(m)}_1,\dots,b^{(m)}_d\big).
\end{gather*}
Since $Q(p)$ is the product of the labels
\begin{gather*}
\prod\limits^d_{i=1}\frac{\tau_{l_{m-s-1}+a^{(s+1)}_i-1}}{\tau_{l_{m-s-1}+a^{(s)}_i}}
\end{gather*}
of the edges $(m-s;a^{(s)}_1,\dots,a^{(s)}_d)\rightarrow (m-s-1;a^{(s+1)}_1,\dots,a^{(s+1)}_d)$, $0\leq s\leq m-1$, the
assumption $Q(p)=Q(p')$ means that
\begin{gather}
\label{abeq}
\prod\limits^d_{i=1}\frac{\tau_{l_{m-s-1}+a^{(s+1)}_i-1}}{\tau_{l_{m-s-1}+a^{(s)}_i}}
=\prod\limits^d_{i=1}\frac{\tau_{l_{m-s-1}+b^{(s+1)}_i-1}}{\tau_{l_{m-s-1}+b^{(s)}_i}}
\end{gather}
for all $0\leq s\leq m-1$.
For $s=0$, by Def\/inition~\ref{pathdef}(iv), we get $a^{(0)}_i=b^{(0)}_i=i$, $1\leq i\leq d$.
It follows from $a^{(1)}_1<a^{(1)}_2<\dots<a^{(1)}_d$, $b^{(1)}_1<b^{(1)}_2<\dots<b^{(1)}_d$ and~\eqref{abeq} for $s=0$
that $a^{(1)}_i=b^{(1)}_i$, $1\leq i\leq d$.
Repeating this argument for $s=1,2,3,\dots,m-1$, we get $a^{(s)}_i=b^{(s)}_i$, $1\leq i\leq d$, which means that $p=p'$.
\end{proof}

\begin{proof}[Proof of Theorem~\ref{thm2}]
We set
\begin{gather*}
\mathbb{B}:=\{Q(p) \,|\, p\in X_d(m,m')\}
=\left\{\prod\limits^{d}_{i=1} \prod\limits_{j=1}^{m-m'}\overline{C}\big(R^{(j)}_i,K^{(j)}_i\big) \,\Bigg|\, R^{(j)}_{i}
:=m-K^{(j)}_{i}-j+i \right\},
\end{gather*}
where $\{K^{(j)}_i\}$ satisf\/ies the conditions
\begin{gather*}
1\leq K^{(j)}_1<K^{(j)}_2<\dots<K^{(j)}_d\leq m'+d,
\qquad
1\leq j\leq m-m',\nonumber
\\
i\leq K^{(1)}_i\leq \dots\leq K^{(m-m')}_i\leq m'+i,
\qquad
1\leq i\leq d.%\label{conditions(*)}
\end{gather*}
By Proposition~\ref{monoprop} and Lemma~\ref{coefflem}, we need to show that
\begin{gather}
\label{lastclaim}
\mathbb{B}=\mu_Y(B^-_{u_{\leq k}}((m'-m)\Lambda_d)),
\end{gather}
where $\mu_Y(x)$ is an embedding of $x\in B^-_{u_{\leq k}}((m'-m)\Lambda_d)$ in Theorem~\ref{thm2}, and
\begin{gather}
\label{yy2}
Y:=\frac{1}{\tau_{l_{m-1}+d}\tau_{l_{m-2}+d}\cdots\tau_{l_{m'}+d}},
\end{gather}
which is the lowest weight vector in $\mu_Y(B^-_{u_{\leq k}}((m'-m)\Lambda_d))$ (Theorem~\ref{monorealmain}(ii)).

First, let us prove the inclusion $\mu_Y(B^-_{u_{\leq k}}((m'-m)\Lambda_d))\subset\mathbb{B}$.
Using $\overline{C}(a,b)$, the monomial~$Y$ in~\eqref{yy2} is described as follows{\samepage
\begin{alignat*}{3}
& &&
\ \  \ \overline{C}(m-1,1)\cdot\overline{C}(m-1,2)\cdots\overline{C}(m-1,d)&  \\
&Y=&&\cdot\ \overline{C}(m-2,1)\cdot\overline{C}(m-2,2)\cdots\overline{C}(m-2,d)&  \\
& & &\qquad \quad \vdots \qquad \qquad \quad \vdots \qquad \qquad \qquad \vdots&\\
& && \cdot\ \overline{C}(m',1)\ \ \cdot \ \ \overline{C}(m',2)\ \ \cdots\ \ \overline{C}(m',d).&
\end{alignat*}

}

\noindent
Thus, we see that $Y\in\mathbb{B}$.
It follows from Theorem~\ref{kashidem} and the def\/inition~\eqref{inc} of $u_{\leq k}$ that
\begin{gather*}
\mu_Y(B^-_{u_{\leq k}}((m'-m)\Lambda_d))=
\big\{\underbrace{\tilde{e}_1^{(N_{1}(m'))}\cdots\tilde{e}_r^{(N_{r}(m'))}}_{\text{$1^{\rm st}$ cycle}}
\underbrace{\tilde{e}_1^{(N_{1}(m'-1))}\cdots\tilde{e}_{r-1}^{(N_{r-1}(m'-1))}}_{\text{$2^{\rm nd}$ cycle}}
\\
\qquad {}\cdots \underbrace{\tilde{e}_1^{(N_1(2))}\cdots\tilde{e}_{r-m+2}^{(N_{r-m+2}(2))}}_{\text{$(m'-1)^{\rm th}$ cycle}}
\underbrace{\tilde{e}_1^{(N_1(1))}\dots\tilde{e}_d^{(N_d(1))}}_{\text{$m'^{\rm th}$ cycle}}\cdot Y \,|\,  N_i(j)\in
\mathbb{Z}_{\geq0}\big\}\setminus\{0\}.
\end{gather*}

By %Lemma~\ref{pathlems}~
\eqref{pathlems2}, for any monomial $Z \in \mathbb{B}$ and $1\leq s < m'+d$, we have
\begin{gather}
\label{esinclude}
\tilde{e}_s Z\in \mathbb{B} \cup \{0\}.
\end{gather}
For an arbitrary set $\{N_i(1)\}_{i=1,\dots,d}$ of non negative integers, the monomial
$\tilde{e}_1^{(N_1(1))}\cdots\tilde{e}_d^{(N_d(1))}\cdot Y$ does not include factors in the form $\overline{C}(*, a)$,
$d+2\leq a\leq r-m+2$, by~\eqref{canda} and %Lemma~\ref{pathlems}
 \eqref{pathlems1}.
Therefore, $\tilde{e}_a\cdot \tilde{e}_1^{(N_1(1))}\cdots\tilde{e}_d^{(N_d(1))}\cdot Y=0$ for $d+2\leq a\leq r-m+2$~by
%Lemma~\ref{pathlems}
\eqref{pathlems1}.
Similarly, for an arbitrary set $\{N_i(2)\}_{i=1,\dots,r-m+2}$ of non negative integers, the monomial
$\tilde{e}_1^{(N_1(2))}\cdots\tilde{e}_{r-m+2}^{(N_{r-m+2}(2))}\tilde{e}_1^{(N_1(1))} $ $\cdots\tilde{e}_d^{(N_d(1))}\cdot Y$
does not include factors in the form $\overline{C}(*, a)$, $d+3\leq a\leq r-m+3$, which means that $\tilde{e}_a\cdot
\tilde{e}_1^{(N_1(2))}\cdots\tilde{e}_{r-m+2}^{(N_{r-m+2}(2))}\tilde{e}_1^{(N_1(1))}$ $\cdots\tilde{e}_d^{(N_d(1))}\cdot Y=0$.

Repeating this argument, we obtain
\begin{gather*}
\mu_Y(B^-_{u_{\leq k}}((m'-m)\Lambda_d))=
\big\{\underbrace{\tilde{e}_1^{(N_{1}(m'))}\cdots\tilde{e}_{m'+d-1}^{(N_{m'+d-1}(m'))}}_{\text{$1^{\rm st}$ cycle}}
\underbrace{\tilde{e}_1^{(N_{1}(m'-1))}\cdots\tilde{e}_{m'+d-2}^{(N_{m'+d-2}(m'-1))}}_{\text{$2^{\rm nd}$ cycle}}
\\
\qquad{}
\cdots \underbrace{\tilde{e}_1^{(N_1(2))}\cdots\tilde{e}_{d+1}^{(N_{d+1}(2))}}_{\text{$(m'-1)^{\rm th}$ cycle}}
\underbrace{\tilde{e}_1^{(N_1(1))}\cdots\tilde{e}_d^{(N_d(1))}}_{\text{$m'^{\rm th}$ cycle}}\cdot Y \,|\,  N_i(j)\in
\mathbb{Z}_{\geq0}\big\}\setminus\{0\}.
\end{gather*}
Therefore, we get $\mu_Y(B^-_{u_{\leq k}}((m'-m)\Lambda_d))\subset\mathbb{B}$ by~\eqref{esinclude}.

Next, we shall prove $\mathbb{B}\subset \mu_Y(B^-_{u_{\leq k}}((m'-m)\Lambda_d))$.
For this, we take an arbitrary element
\begin{gather*}
M:=\prod\limits^{d}_{i=1} \prod\limits_{j=1}^{m-m'}\overline{C}(R^{(j)}_i,K^{(j)}_i)\in \mathbb{B}.
\end{gather*}
We need to show that there exists a~set $\{N_i(s)\}$ of non negative integers such that
\begin{gather}
\underbrace{\tilde{e}_1^{(N_{1}(m'))}\cdots\tilde{e}_{m'+d-1}^{(N_{m'+d-1}(m'))}}_{\text{$1^{\rm st}$ cycle}}
\underbrace{\tilde{e}_1^{(N_{1}(m'-1))}\cdots\tilde{e}_{m'+d-2}^{(N_{m'+d-2}(m'-1))}}_{\text{$2^{\rm nd}$ cycle}}
\nonumber\\
\qquad {}\cdots \underbrace{\tilde{e}_1^{(N_1(2))}\cdots\tilde{e}_{d+1}^{(N_{d+1}(2))}}_{\text{$(m'-1)^{\rm th}$ cycle}}
\underbrace{\tilde{e}_1^{(N_1(1))}\cdots\tilde{e}_d^{(N_d(1))}}_{\text{$m'^{\rm th}$ cycle}}\cdot Y=M.  \label{lastclaim2}
\end{gather}

We set $\{N_i(s)\}$ as follows (see Example~\ref{pathex4}). First, we set $\{N_i(1)\}_{i=1,2,\dots,d}$ as
\begin{gather*}
N_d(1) := \#\big\{1\leq j\leq m-m' \,|\,  K^{(j)}_d -d =m' \big\},
\\
N_{d-1}(1) := \#\big\{1\leq j\leq m-m' \,|\,  K^{(j)}_{d-1}-(d-1)=m' \big\},
\\
N_{d-2}(1) := \#\big\{1\leq j\leq m-m' \,|\,  K^{(j)}_{d-2}-(d-2)=m' \big\},
\\
 \cdots\cdots\cdots\cdots\cdots\cdots\cdots\cdots\cdots\cdots\cdots\cdots\cdots\cdots\cdots
\\
N_{1}(1) := \#\big\{1\leq j\leq m-m' \,|\,  K^{(j)}_{1}-1=m' \big\}.
\end{gather*}
Note that since $K^{(j)}_1<\dots<K^{(j)}_d$, we have $N_1(1)<\dots<N_d(1)$.
As seen in Example~\ref{exex}, by applying $\tilde{e}^{(N_d(1))}_d$ to~$Y$, the factors $\overline{C}(m',d)$,
$\overline{C}(m'+1,d),\dots,\overline{C}(m'+N_d(1)-1,d)$ in~$Y$ turn out $\overline{C}(m'-1,d+1)$,
$\overline{C}(m',d+1),\dots,\overline{C}(m'+N_d(1)-2,d+1)$ by Remark~\ref{actionrem}, \eqref{canda} and
%Lemma~\ref{pathlems}
\eqref{pathlems1}.
Then the monomial $\tilde{e}^{(N_d(1))}_d\cdot Y$ does not include factors in the form $\tau_{l_q+d-1}$, $0\leq q\leq m-1$.
Thus, we can use Example~\ref{exex} again when we apply $\tilde{e}^{(N_{d-1}(1))}_{d-1}$ to
$(\tilde{e}^{(N_d(1))}_d\cdot Y)$.
Then the factors $\overline{C}(m',d-1)$, $\overline{C}(m'+1,d-1),\dots,\overline{C}(m'+N_{d-1}(1)-1,d-1)$ in
$(\tilde{e}^{(N_d(1))}_d\cdot Y)$ turn out $\overline{C}(m'-1,d)$,
$\overline{C}(m',d),\dots,\overline{C}(m'+N_{d-1}(1)-2,d)$, and the monomial
$\tilde{e}^{(N_{d-1}(1))}_{d-1}\tilde{e}^{(N_d(1))}_d\cdot Y$ does not include factors in the form $\tau_{l_q+d-2}$, $0\leq q\leq m-1$.
After all, by using Example~\ref{exex} repeatedly, we see that the monomial
$\tilde{e}^{(N_{1}(1))}_{1}\cdots\tilde{e}^{(N_{d-1}(1))}_{d-1}\tilde{e}^{(N_d(1))}_d\cdot Y$ is obtained from~$Y$~by
replacing $\overline{C}(m',\zeta)$, $\overline{C}(m'+1,\zeta),\dots,\overline{C}(m'+N_{\zeta}(1)-1,\zeta)$~by
$\overline{C}(m'-1,\zeta+1)$, $\overline{C}(m',\zeta+1),\dots,\overline{C}(m'+N_{\zeta}(1)-2,\zeta+1)$, $1\leq\zeta\leq d$.

We set $\{N_i(s)\}_{2\leq s \leq m',\ i=1,2,\dots,d+s-1}$ as
\begin{gather*}
N_{d+1}(2) := N_d(1),
\\
N_d(2) := N_{d-1}(1)+\#\big\{1\leq j\leq m-m' \,|\,  K^{(j)}_d -d =m'-1 \big\},
\\
N_{d-1}(2) := N_{d-2}(1)+\#\big\{1\leq j\leq m-m' \,|\,  K^{(j)}_{d-1}-(d-1)=m'-1 \big\},
\\
 \cdots\cdots\cdots\cdots\cdots\cdots\cdots\cdots\cdots\cdots\cdots\cdots\cdots
\\
N_{2}(2) := N_{1}(1)+\#\big\{1\leq j\leq m-m' \,|\,  K^{(j)}_{2}-2=m'-1 \big\},
\\
N_{1}(2) := \#\big\{1\leq j\leq m-m' \,|\,  K^{(j)}_{1}-1=m'-1 \big\},
\\
N_{d+2}(3) := N_{d+1}(2),
\\
N_{d+1}(3) := N_{d}(2),
\\
N_d(3) := N_{d-1}(2)+\#\big\{1\leq j\leq m-m' \,|\,  K^{(j)}_d -d =m'-2 \big\},
\\
N_{d-1}(3) := N_{d-2}(2)+\#\big\{1\leq j\leq m-m' \,|\,  K^{(j)}_{d-1}-(d-1)=m'-2 \big\},
\\
 \cdots\cdots\cdots\cdots\cdots\cdots\cdots\cdots\cdots\cdots\cdots\cdots\cdots
\\
N_{2}(3) := N_{1}(2)+\#\big\{1\leq j\leq m-m' \,|\,  K^{(j)}_{2}-2=m'-2 \big\},
\\
N_{1}(3) := \#\big\{1\leq j\leq m-m' \,|\,  K^{(j)}_{1}-1=m'-2 \big\},
\\
 \cdots\cdots\cdots\cdots\cdots\cdots\cdots\cdots\cdots\cdots\cdots\cdots\cdots
\\
N_{m'+d-1}(m') := N_{m'+d-2}(m'-1),
\\
N_{m'+d-2}(m') := N_{m'+d-3}(m'-1),
\\
 \cdots\cdots\cdots\cdots\cdots\cdots\cdots\cdots\cdots\cdots\cdots\cdots\cdots
\\
N_{d+1}(m') := N_d(m'-1),
\\
N_d(m') := N_{d-1}(m'-1)+\#\big\{1\leq j\leq m-m' \,|\,  K^{(j)}_d-d =1 \big\},
\\
 \cdots\cdots\cdots\cdots\cdots\cdots\cdots\cdots\cdots\cdots\cdots\cdots\cdots
\\
N_2(m') := N_{1}(m'-1)+\#\big\{1\leq j\leq m-m' \,|\,  K^{(j)}_2-2 =1 \big\},
\\
N_1(m') := \#\big\{1\leq j\leq m-m' \,|\,  K^{(j)}_1-1=1 \big\}.
\end{gather*}
For example, if $K^{(m-m')}_d-d=m'$, then the factor $\overline{C}(m',d)$ in~$Y$ is acted~by
$\tilde{e}_{m'+d-1}\cdots\tilde{e}_{d+2}\tilde{e}_{d+1}\tilde{e}_d$ since
$N_{m'+d-1}(m')=N_{m'+d-2}(m'-1)=N_{m'+d-3}(m'-2)=\cdots=N_{d}(1)$, and we obtain
\begin{gather*}
\tilde{e}_{m'+d-1}\cdots\tilde{e}_{d+2}\tilde{e}_{d+1}\tilde{e}_d\cdot
\overline{C}(m',d)=\overline{C}(0,m'+d)=\overline{C}\big(R^{(m-m')}_d,K^{(m-m')}_d\big),
\end{gather*}
by~\eqref{canda} and %Lemma~\ref{pathlems}
\eqref{pathlems1}.
In general, if $K^{(j)}_i-i=m'-\zeta$, $0\leq\zeta\leq m'-1$, then the factor $\overline{C}(m-j,i)$ in~$Y$ is acted~by
$\tilde{e}_{m'-\zeta+i-1}\cdots\tilde{e}_{i+2}\tilde{e}_{i+1}\tilde{e}_i$, and we obtain
\begin{gather*}
\tilde{e}_{m'-\zeta+i-1}\cdots\tilde{e}_{i+2}\tilde{e}_{i+1}\tilde{e}_i\overline{C}(m-j,i)
=\overline{C}(m-m'-j+\zeta,m'-\zeta+i)=\overline{C}(R^{(j)}_i,K^{(j)}_i),
\end{gather*}
which means~\eqref{lastclaim2}.
Therefore, we get~\eqref{lastclaim}.
\end{proof}

\begin{ex}
\label{pathex4}
We set $r=9$, $m=6$, $m'=3$ and $d=4$.
Let us see that we can obtain
\begin{alignat*}{3}
 && & \overline{C}(6,1)\cdot\overline{C}(6,2)\cdot \overline{C}(5,4)\cdot \overline{C}(4,6)&
\\
& M:={} && \overline{C}(4,2)\cdot\overline{C}(4,3)\cdot \overline{C}(4,4)\cdot \overline{C}(3,6)&
\\
 &&& \overline{C}(1,4)\cdot\overline{C}(1,5)\cdot \overline{C}(1,6)\cdot \overline{C}(1,7)  &
\end{alignat*}
from
\begin{alignat*}{3}
&&&  \overline{C}(6,1)\cdot\overline{C}(6,2)\cdot \overline{C}(6,3)\cdot \overline{C}(6,4)&
\\
 &Y:={} && \overline{C}(5,1)\cdot\overline{C}(5,2)\cdot \overline{C}(5,3)\cdot \overline{C}(5,4)&
\\
&&&  \overline{C}(4,1)\cdot\overline{C}(4,2)\cdot \overline{C}(4,3)\cdot \overline{C}(4,4).&
\end{alignat*}
by applying $\{\tilde{e}_i^{N_i(s)}\}$ in the proof of Theorem~\ref{thm2}, that is
\begin{gather}
\underbrace{\tilde{e}_{1}^{(N_1(3))}\tilde{e}_{2}^{(N_2(3))}\tilde{e}_{3}^{(N_3(3))}\tilde{e}_{4}^{(N_4(3))}
\tilde{e}_{5}^{(N_5(3))}\tilde{e}_{6}^{(N_6(3))}}_{\text{$3^{\rm rd}$ cycle}} \cdot
\underbrace{\tilde{e}_{1}^{(N_1(2))}\tilde{e}_{2}^{(N_2(2))}\tilde{e}_{3}^{(N_3(2))}\tilde{e}_{4}^{(N_4(2))}
\tilde{e}_{5}^{(N_5(2))}}_{\text{$2^{\rm nd}$ cycle}} \nonumber
\\
\qquad {}\cdot \underbrace{\tilde{e}_{1}^{(N_1(1))}\tilde{e}_{2}^{(N_2(1))}\tilde{e}_{3}^{(N_3(1))}\tilde{e}_{4}^{(N_4(1))}}_{\text{$1^{\rm st}$ cycle}}\cdot Y=M.\label{ymeq}
\end{gather}

We set $M=\prod\limits^3_{j=1}\prod\limits^4_{i=1}\overline{C}(R^{(j)}_i,K^{(j)}_i)$, that is, $K^{(1)}_1=1$,
$K^{(1)}_2=2$, $K^{(1)}_3=4$, $K^{(1)}_4=6$, $K^{(2)}_1=2$, $K^{(2)}_2=3, \dots, K^{(3)}_4=7$.

To change the factor $\overline{C}(4,4)$ of~$Y$ into $\overline{C}(1,7)$ of~$M$, we need to apply
$\tilde{e}_6\tilde{e}_5\tilde{e}_4$ to $\overline{C}(4,4)$.
In contrast, to change the factors $\overline{C}(5,4)$ and $\overline{C}(6,4)$ of~$Y$ into $\overline{C}(3,6)$ and
$\overline{C}(4,6)$, we need to apply $\tilde{e}_5\tilde{e}_4$ to $\overline{C}(5,4)$ and $\overline{C}(6,4)$.
Thus, $\overline{C}(4,4)$ must be changed to $\overline{C}(3,5)$ by the action of~$\tilde{e}_{4}^{(N_4(1))}$ in f\/irst
cycle, and $\overline{C}(5,4)$, $\overline{C}(6,4)$ do not have to be changed at the f\/irst cycle.
So we set $N_4(1):=1=\#\{j \,|\, K^{(j)}_4-4=3\}$.
Similarly, we set $N_3(1)=N_2(1)=N_1(1)=1$, and the factors $\overline{C}(4,1)$, $\overline{C}(4,2)$, $\overline{C}(4,3)$
of~$Y$ is changed to $\overline{C}(3,2)$, $\overline{C}(3,3)$, $\overline{C}(3,4)$ by the action of f\/irst cycle.

Next, to obtain~$M$, the factor $\overline{C}(3,5)$ in
$\tilde{e}_{1}^{(N_1(1))}\tilde{e}_{2}^{(N_2(1))}\tilde{e}_{3}^{(N_3(1))}\tilde{e}_{4}^{(N_4(1))}\cdot Y$ must be
changed to $\overline{C}(2,6)$ at the second cycle.
Thus, we set $N_5(2):=1=N_4(1)$.
The factor $\overline{C}(3,4)$ must be changed to $\overline{C}(2,5)$ by the action of $\tilde{e}_{4}^{(N_4(2))}$ in
second cycle, and the factors $\overline{C}(5,4)$ and $\overline{C}(6,4)$ must also be changed at the second cycle.
Hence, we set $N_4(2):=3=N_3(1)+\#\{j \,|\, K^{(j)}_4-4=2\}$.
Similarly, we set $N_3(2)=1$, $N_2(2)=1$, $N_1(2)=0$, $N_6(3)=1$, $N_5(3)=3$, $N_4(3)=1$, $N_3(3)=3$, $N_2(3)=1$ and
$N_1(3)=1$.
Then we get~\eqref{ymeq}.
\end{ex}

\begin{ex}
%\label{pathex5}
We use the same setting in Example~\ref{pathex3}: $u=s_1s_2s_3s_4s_1s_2s_3s_1s_2s_1$, $k=6$, $m=3$ and $m'=2$.
We have $u_{\leq6}=s_1s_2s_3s_4s_1s_2$.
Let
\begin{gather*}
Y:=\frac{1}{\tau_{l_{2}+2}}\in{\mathcal Y},
\end{gather*}
which has weight $-\Lambda_2$.
By Theorem~\ref{thm2}, we obtain
\begin{gather*}
\Delta^L(6;{\bf i})(\tau) =  \sum\limits_{x\in B^-_{u_{\leq 6}}(-\Lambda_2)}\mu_Y(x)
 =  Y+\tilde{e}_2 Y+\tilde{e}_1\tilde{e}_2 Y+\tilde{e}_3\tilde{e}_2 Y +\tilde{e}_3\tilde{e}_1\tilde{e}_2 Y
+\tilde{e}_2\tilde{e}_3\tilde{e}_1\tilde{e}_2 Y
\\
\phantom{\Delta^L(6;{\bf i})(\tau)}
 = \frac{1}{\tau_{l_{2}+2}}+\frac{\tau_{l_1+2}}{\tau_{l_2+1}\tau_{l_1+3}}
+\frac{\tau_{l_1+1}}{\tau_{l_1+3}}+\frac{\tau_{l_0+3}}{\tau_{l_2+1}\tau_{l_0+4}}
+\frac{\tau_{l_1+1}\tau_{l_0+3}}{\tau_{l_1+2}\tau_{l_0+4}} +\frac{\tau_{l_0+2}}{\tau_{l_0+4}}.
\end{gather*}
\end{ex}

\subsection{The proof of Theorem~\ref{thm3}}\label{prthm3}

Let us prove Theorem~\ref{thm3}.
Suppose that $i_k=d=1$.

\begin{proof}[Proof of Theorem~\ref{thm3}]
By Proposition~\ref{pathlem}, we have
\begin{gather}
\label{ik1}
\Delta^L(k;{\bf i})(\tau)=\sum\limits_{p\in X_1(m,m')} Q(p).
\end{gather}

The set $X_1(m,m')$ consists of paths~$p$ as in the form
\begin{gather*}
p=(m,1)\rightarrow\big(m-1,a^{(1)}\big)\rightarrow\big(m-2,a^{(2)}\big)\rightarrow\dots \rightarrow\big(1,a^{(m-1)}\big)\rightarrow(0,m'+1)
\end{gather*}
such that $a^{(s+1)}=a^{(s)}$ or $a^{(s)}+1$, $0\leq s\leq m-1$.
Here, $a^{(0)}:=1$, $a^{(m)}:=m'+1$.

By Lemma~\ref{nmblem}, we obtain
\begin{gather*}
\#\big\{s \,|\, a^{(s+1)}=a^{(s)}+1,\ 0\leq s\leq m-1 \big\}=m'.
\end{gather*}
Set $\{s \,|\, a^{(s+1)}=a^{(s)}+1,\; 0\leq s\leq m-1 \}:=\{j_1,\dots,j_{m'}\}$, $0\leq j_1<\dots<j_{m'}\leq m-1$.
Then we have
\begin{gather}
a^{(0)}=1,\quad  a^{(1)}=1,\quad  a^{(2)}=1,\quad \dots,\quad a^{(j_1)}=1,\quad  a^{(j_1+1)}=2,
\nonumber
\\
a^{(j_1+2)}=2,\quad  a^{(j_1+3)}=2,\quad \dots,\quad  a^{(j_2)}=2,\quad  a^{(j_2+1)}=3,
\nonumber
\\
a^{(j_2+2)}=3,\quad  a^{(j_2+3)}=3,\quad \dots, \quad a^{(j_3)}=3,\quad  a^{(j_3+1)}=4,
\nonumber
\\
\label{jlist4}
\cdots\cdots\cdots\cdots\cdots\cdots\cdots\cdots\cdots\cdots\cdots\cdots\cdots\cdots\cdots\cdots\cdots\cdots
\\
\nonumber
a^{(j_{\nu}+2)}=\nu+1,\quad  a^{(j_{\nu}+3)}=\nu+1,\quad \dots, \quad a^{(j_{\nu+1})}=\nu+1,\quad  a^{(j_{\nu+1}+1)}=\nu+2,
\\
\cdots\cdots\cdots\cdots\cdots\cdots\cdots\cdots\cdots\cdots\cdots\cdots\cdots\cdots\cdots\cdots\cdots\cdots
\nonumber
\\
a^{(j_{m'}+2)}=m'+1,\quad  a^{(j_{m'}+3)}=m'+1,\quad \dots, \quad a^{(m)}=m'+1.
\nonumber
\end{gather}
Therefore, Def\/inition~\ref{labeldef}(ii) means
\begin{gather*}
Q(p)=\prod\limits^{j_1-1}_{i=0}\frac{1}{\tau_{l_{m-1-i}+1}}
\prod\limits^{j_2-1}_{i=j_1+1}\frac{\tau_{l_{m-1-i}+1}}{\tau_{l_{m-1-i}+2}}\cdots
\prod\limits^{m-1}_{i=j_{m'}+1}\frac{\tau_{l_{m-1-i}+m'}}{\tau_{l_{m-1-i}+m'+1}}.
\end{gather*}

Conversely, for a~given $\{j_1,\dots,j_{m'}\}$, $0\leq j_1<\dots<j_{m'}\leq m-1$, we can constitute a~path~$p$ as in~\eqref{jlist4}.
Hence, by~\eqref{ik1}, we proved our claim.
\end{proof}

\subsection*{Acknowledgements}
The authors would like to acknowledge the referees for giving them relevant advice and suggestion to improve this
article.
T.N.~is supported in part by JSPS Grants in Aid for Scientif\/ic Research $\sharp 22540031$, $\sharp 15K04794$.

\pdfbookmark[1]{References}{ref}
\LastPageEnding

\end{document}